\documentclass[11pt]{article}
\usepackage{amsfonts}
\usepackage{mathrsfs}
\usepackage{amssymb}
\usepackage{amsmath}
\usepackage{algorithm}
\usepackage{cite}
\usepackage{color}
\usepackage{multirow}
\usepackage{extarrows}
\usepackage{graphicx}
\usepackage{subfigure}
\usepackage{float}
\usepackage{diagbox}
\usepackage{booktabs}
\usepackage{algpseudocode}
\usepackage{caption}
\usepackage{comment}
\usepackage{adjustbox}
\usepackage{graphicx}

\allowdisplaybreaks[4]

\newtheorem{theorem}{Theorem}

\newcommand\btd{\raise 2pt \hbox{$\hat\bigtriangledown$}\hskip 1.5pt}
\newcommand\bt{\raise 2pt \hbox{$\bigtriangledown$}\hskip 1.5pt}

\def\no{\nonumber}
\def\x{\textbf{x}}
\def\b{\textbf{b}}

\def\w{\textbf{w}}

\begin{document}
\title{HANN: Homotopy auxiliary neural network for solving nonlinear algebraic equations}
\author{Ling-Zhe Zai $^1$ \ \ \  Lei-Lei  Guo $^2$ \footnote{Ling-Zhe Zai and Lei-Lei  Guo are co-first authors of the article.} \ \ \  Zhi-Yong Zhang $^1$\footnote{E-mail: zzy@muc.edu.cn}
 \\
\small $^1$ College of Science, Minzu University of China, Beijing 100081, P.R. China\\
\small $^2$ College of Science, North China University of Technology, Beijing 100144, P.R. China}
\date{}
\maketitle
\noindent{\bf Abstract:} Solving nonlinear algebraic equations is a fundamental but challenging problem in scientific computations and also has many applications in system engineering. Though traditional iterative methods and modern optimization algorithms have exerted effective roles in addressing certain specific problems, there still exist certain weaknesses such as the initial value sensitivity, limited accuracy and slow convergence rate, particulary without flexible input for the neural network methods. In this paper, we propose a homotopy auxiliary neural network (HANN) for solving nonlinear algebraic equations which integrates the classical homotopy continuation method and popular physics-informed neural network. Consequently, the HANN-1 has strong learning ability and can rapidly give an acceptable solution for the problem which outperforms some known methods, while the HANN-2 can further improve its accuracy. Numerical results on the benchmark problems confirm that the HANN method can effectively solve the problems of determining the total number of solutions of a single equation, finding solutions of transcendental systems involving the absolute value function or trigonometric function, ill-conditioned and normal high-dimensional nonlinear systems and time-varying nonlinear problems, for which the Python's built-in Fsolve function exhibits significant limitations, even fails to work. 



\section{Introduction}

Nonlinear algebraic equations arise in many fields of practical importance such as engineering, mechanics, chemistry and robotics. Thus finding solutions of a system of nonlinear algebraic equations is an important and  fundamental topic in science and contributes to the various applied subjects, for example, geometric intersection problems \cite{dd-1986}, computation of equilibrium states \cite{bh-2010}, etc., but it is a tough task and up to now no unified and efficient method exists to tackle it. Currently, there are two primary approaches of symbolic  computation and numerical computation. From the point of view of symbolic computation, the elimination of variables is the usual adopted approach, such as Gr\"{o}bner bases \cite{bb-1985}, characteristic set \cite{wu}, resultant theory \cite{cox}. However, such symbolic computation methods suffer from high computational complexities because of performing variable elimination and thus are only limited to low-degree systems of polynomials, and fail to handle a little complicated algebraic equations in practical computational environments. For example, when computing the Gr\"{o}bner bases of a polynomial system with three variables, there exists an intermediate polynomial whose four coefficients contain roughly eighty thousand digits \cite{ar-2003}.

Therefore, numerical methods based on iterative idea dominate the overwhelming superiority in solving nonlinear algebraic equations. For example, the classical Newton-Raphson method has the quadratic convergence rate near roots, but requires computing the Jacobian matrix and is also sensitive to the initial guess \cite{kell}. The homotopy continuation method, originated from Garcia and Zangwill \cite{gz-1979} and Drexler \cite{dre-1977} independently, starts with a trivial system of easy solving and gradually evolves into the target nonlinear system to be solved by moving on a path of the constructed continuous homotopy. The main advantage of homotopy continuation method is that it can solve the large degree of polynomial system and has the ability to find all isolated solutions of the desired system, even including singular and multiple roots. However, it strongly depends on the initial solutions of trivial system and also suffers from the path tracking problem where the solution path may diverge or cross, as well as the high computational cost for the high-dimensional or highly nonlinear systems \cite{li-1997}. Alternatively, an evolutionary computation technique was used to solve the nonlinear equations by converting the system into a multi-objective optimization problem \cite{gaa-2008}, whereas the multi-objective optimization itself is a hard problem to be solved. In addition, the clustering-based methods \cite{wsn-2011} and repulsion-based methods \cite{gwz-2020} were proposed to find the multiple solutions of nonlinear algebraic equations. However, such methods still have their limitations though they have proven effective in the specific problems. For example, the clustering methods require manual specification of cluster numbers or rely on error-prone heuristics, while the multi-objective optimization methods are plagued by the curse of dimensionality and pose significant challenges in the optimization field.

Another main method for solving the algebraic equations numerically is the neural network based methods which work on the foundation of the universal approximation theory \cite{hsm-1989}, and  transform the studied equations into an optimization problem where the solutions are approximated by training the weight matrixes and bias vectors of neural network. Such neural network methods contain two main steps of constructing network model by introducing the loss function and training the model with gradient descent, and have some obvious merits such as without needing explicit mathematical expressions, handling high-dimensional nonlinear systems and performing with parallel computations \cite{tt-1993,mk-2007,ka-2017}. For example, Nguyen implemented the well-known Newton-Raphson algorithm on the neural network architecture where the computational speed is very fast for solving large nonlinear system and the computational time is independent of the dimension of system \cite{tt-1993}. Mishra and Kalra proposed a modified Hopfield neural network where a new energy function was introduced and used to update the weight matrixes and bias vectors \cite{mk-2007}. A neural network framework with an adaptive learning rate for solving the polynomial equations was proposed in \cite{ka-2017}. 
However, the current neural network based methods for solving algebraic equations typically reformulate the iterative methods in the framework of neural network because the algebraic equations themselves only involve the unknown variables and thus force the methods to use the same variables as both inputs and outputs. Such a design lacks auxiliary input variables and thereby limits the framework's adaptability.


In fact, when solving differential equations with the neural network method, the independent variables are chosen as the inputs and the dependent variables serve as the outputs, which grants the network greater flexibility during training, thereby enhancing its learning ability \cite{2018a,gm-2018}. Among these methods, physics-informed neural networks (PINN) stands out by integrating residuals of both the differential equations and initial/boundary conditions into the loss function \cite{2018a}. Up to now, the PINN and its improved versions, for example, incorporating certain additional physical constraints into the network architecture or the loss function of PINN \cite{zhang-2025a,li-2024,zhang-2022,tb-2022}, have exerted significant effects in fluid mechanics and many real world applications \cite{lmz-2021,mag-2020}. Therefore, a key point for solving algebraic equations is how to construct a neural network framework where the input and output variables are different, analogous to the differential equation case.

In this paper, we present a homotopy auxiliary neural network (HANN) framework for solving algebraic equations, which first incorporates the target algebraic system into a homotopy model, then borrows the idea of PINN for solving differential equations to train the homotopy model where the introduced homotopy parameter serves as an input variable. More specifically, we introduce two HANN methods, denoted by HANN-1 and HANN-2, whose main advantages include:

$\bullet$ The HANN-1 combines the homotopy model connecting the original system with an easy solving system and an initial condition into the loss function of neural network and operates purely through a standard feed-forward propagation, removing dependencies on iterative solvers or manual initialization tuning.

$\bullet$ Compared with the conventional neural network-based methods for algebraic equations where the input is the user-specified initial values of unknown variables, the two HANN methods employ the homotopy parameter as the input variable and thus offer flexibility on the properties of input, such as the number of input data and the sampling density of defined domain.

$\bullet$ The HANN-1 exhibits strong learning capability to solve nonlinear algebraic equations including the high-dimensional polynomial equations, the transcendental equations (featuring non-algebraic functions like trigonometric or exponential terms) and the time-varying equations, for which the Python's built-in Fsolve function fails. Moreover, the computational costs of HANN-1 is highly low.

$\bullet$ The HANN-2 further improves the accuracies of the learned solutions and outperforms the HANN-1 and some state-of-the-art methods.

The remainder of the paper is organized as follows. In Section 2, we introduce the main idea of HANN where the theories and algorithms of HANN-1 and HANN-2 are presented in detail. In Section 3, we illustrate the efficacies of the two HANN methods by performing numerical experiments for the several benchmark examples, particularly for the time-varying nonlinear equations. The last section concludes the results.

\section{Methods}
By means of homotopy continuation method and physics-informed neural network, we propose a novel framework of deep neural network to find solutions of a nonlinear algebraic system
\begin{equation}\label{aleq}
F(\x):=(f_1(\x),f_2(\x),\dots,f_n(\x))=0,
\end{equation}
where $\x=(x_1,x_2,\dots,x_n) \in \prod_{i=1}^n[\alpha_i,\beta_i]\subset \mathbb{R}^n$ is a vector of unknowns, and the nonlinear system means that at least one equation in $F(\x)=0$ is nonlinear. 

Specifically, by introducing the parameter $t$, we construct a homotopy model related with system (\ref{aleq}) by
\begin{equation}\label{h-model}
H(\x,t):= t F(\x) + \gamma  (t-1)  (F(\x)-F(\x_0)),
\end{equation}
where $\x_0$ is an initial value of $\x$ and $\gamma$ is a convergence-control parameter to be adjusted. In particular, 
\begin{eqnarray}
&&\no t=0:~~~H(\x,0) = -\gamma (F(\x)- F(\x_0));\\
&&\no t=1:~~~H(\x,1) = F(\x).
\end{eqnarray}

Then, the starting system $H(\x_0,0)=0$ holds identically for any given initial value $\x_0$ and thus any initial value $\x_0$ which keeps $F(\x)$ nonsingular can be chosen as the starting point. Meanwhile, solving $F(\x)=0$ is equivalent to first find solution $\x(t)$ of the final system $H(\x,t)=0$ and then by setting $t=1$ to obtain the required solution $\x(1)$. Furthermore, Theorem \ref{th-im} assures that we can obtain solution of system (\ref{aleq}) through solving the homotopy system (\ref{h-model}). 
\begin{theorem}[\cite{1}]
Let $U$ be an open subset of $\mathbb{R}^n$ and $\x\in U$. Suppose that $H(\x,t):=(H_1(\x,t),$ $H_2(\x,t),\dots,H_n(\x,t))$ is continuously differentiable in an open set containing $U \times[0,1]$, that the function $ H(\x,0)$ has a zero $\x_0$ in $U$, and that a constant $d > 0$ makes $\{\x: |\x - \x_0| < d\} \subseteq U$. Moreover, on $ U\times[0, 1]$, define two matrixes
\[
D_x H = \left( \frac{\partial H_i}{\partial \x_j} \right)_{i, j = 1, 2, \ldots, n},~~~~D_t H = \left( \frac{\partial H_i}{\partial t} \right)_{i = 1, 2, \ldots, n}.
\]
If $D_x H$ is nonsingular and $\left| (D_x H)^{-1} D_t H \right| < d$ holds, then there exists a continuously differentiable function $\x(t)$ satisfying
$H(\x(t),t)=0$ for $0 \leq t \leq 1$. In particular, we have $F(\x(1)) = H(\x(1),1) = 0.$
\label{th-im}
\end{theorem}

Next, we introduce a neural network framework to solve system (\ref{aleq}) based on the homotopy system (\ref{h-model}) and call it by HANN-1. Since the solution $\x$ of $H(\x,t)=0$ is a function of $t$,  we feed the neural network with the data of $t\in[0,1]$ and take $\x$ as the output. Then the proposed HANN-1 is a deep neural network $\widehat{\x}(\Theta; t)$ and aims to approximate the solution of system (\ref{aleq}), where \(\Theta = \{\mathbf{w}_k, \mathbf{b}_k\}_{k=1}^K\) is the collection of trainable parameters including weight matrices \(\mathbf{w}_k\) and bias vectors \(\mathbf{b}_k\). Specifically, for a HANN-1 with $K$ layers and $n_k$ neurons in the $k$th layer, then for the input data $t_0$, the output of HANN-1 is
\begin{align}
\no \widehat{x}(\Theta; t_0) = \mathcal{F}_K \circ \sigma \circ \mathcal{F}_{K-1} \circ \cdots \circ \sigma \circ \mathcal{F}_1 (t_0),
\end{align}
where the output of the $k$-th layer $t_k = \sigma\left(\mathcal{F}_k(t_{k-1})\right) = \sigma\left(\mathbf{w}_k t_{k-1} + \mathbf{b}_k\right)$, $k=1,\dots,K$,  and $\sigma(\cdot)$ denotes the activation function.

In order to optimize the training parameters $\Theta$, the HANN-1 utilizes the Adam or L-BFGS optimization algorithms \cite{km-2015,ln-1989} to minimize the loss function defined by
\begin{eqnarray}\label{loss}
\mathcal{L}(\Theta) = \varpi_i \mathcal{L}_{iv}(\Theta) +  \varpi_h \mathcal{L}_h(\Theta),
\end{eqnarray}
where \(\varpi_i\) and \(\varpi_h\) are the weights of the loss terms and set to one uniformly, \(\mathcal{L}_{iv}\) measures the degree of data-fit for the unique initial data \(\{(0,\x_0)\}\), while \(\mathcal{L}_h\) is a loss term that penalizes system (\ref{h-model}) on a finite set of collocation points \(\{ t_j\}_{j=1}^{N_f}\) that are selected from the interior domain by the Latin hypercube sampling \cite{ms-1987}. Specifically,
\begin{eqnarray} \label{eq:MSE_f}
\no &&\mathcal{L}_{iv}(\Theta)= \left| \hat{\x}(\Theta; 0) - \x_0 \right|^2, ~~~\mathcal{L}_h(\Theta)= \frac{1}{N_f} \sum_{j=1}^{N_f}  \sum_{i=1}^{n}\left| h_i(\Theta;\hat{\x}_j,t_j) \right|^2,
\end{eqnarray}
where $H(\x,t)=(h_1(\x,t),h_2(\x,t),\dots,h_n(\x,t))$. For clarity, we state the training steps of HANN-1 in Algorithm \ref{alg:HANN_framework}. Furthermore, the HANN-1 formulates the problem of solving system (\ref{aleq}) into a new framework of neural network and is implemented in a standard feed-forward mode, thus eliminates the dilemma that the conventional neural network methods for solving system (\ref{aleq}) only employ the same input and output variables.
\begin{algorithm}[ht]
\caption{Framework of HANN-1}
\label{alg:HANN_framework}
\begin{algorithmic}[1]
    \Require
    Initial value $\x_0$; Equations $F(\x)=0$ and its domain; Number of collocation points \( N_f \);
     \Ensure
      Predicted solution $\x(1)$; Equation residuals of $F(\x(1))$;
    \State Construct the homotopy model (\ref{h-model}) associated with $F(\x)=0$;
    \State Initialize the weight matrices $\w_k$ and bias vectors $\b_k$ $(k=1,\dots,K)$ using the Xavier initialization method;
    \State Construct the loss function $\mathcal{L}(\Theta)$ given by (\ref{loss}) as the sum of mean squared errors of the initial conditions and the residuals of $F(\x)$;
    \State Train the neural network using the Adam optimizer or L-BFGS optimizer;
    \State \Return  Predicted solution $\x(1)$; Equation residuals of $F(\x(1))$;
\end{algorithmic}
\end{algorithm}

Furthermore, it is observed that the HANN-1 can give the learned solutions closer to the true solution than the initial input. Moreover, Theorem \ref{th-im} implies that when the initial value lies in the neighbourhood with a small radius $d$, it is relatively easy to approximate the true solution. Therefore, both the theory in Theorem  \ref{th-im} and the algorithm of HANN-1 motivate us to combine the iteration idea and the HANN-1 to propose HANN-2 for obtaining high-accuracy solutions of system (\ref{aleq}), whose explicit steps are stated in Algorithm \ref{alg:HANN-2}. Then the HANN-2 operates by repeatedly applying HANN-1 with a fixed number of iterations, and terminates if the equation residual increases for ten consecutive steps or the maximum iteration count reaches the given number, which ensures that the HANN-2 consistently outperforms HANN-1 in terms of solution accuracy.

\begin{algorithm}[ht]
\caption{Framework of  HANN-2}
\label{alg:HANN-2}
\begin{algorithmic}[100]
    \Require
        Initial conditions $\mathbf{x}_0$; Equations $F(\mathbf{x}) = 0$ and its domain; Maximal iteration number $N_m$; Number of collocation points $N_f$;
    \Ensure
       Predicted solution \( ans_\mathbf{x} \); Equation residuals of \(ans\_f\);
        \State Let $ans\_f = 1000$; $num_t = num_l =0$;
    \While{\( num_t \leq N_m\) }
        \State Call HANN-1 and train it with the initial conditions $\x_0$, then output $\mathbf{x}(1)$ and the equation residuals $ans\_F$ of $F(\mathbf{x})$
        \If{\( ans\_F \leq ans\_f \)}
            \State \( ans_\mathbf{x} \) = $\mathbf{x}$(1); $\mathbf{x}_0$ = $\mathbf{x}$(1); $ans\_f = ans\_F$
           \State $num_l\leftarrow 0$
        \Else
         \State $num_l\leftarrow num_l + 1$
         \EndIf
        \State $num_t \leftarrow num_t + 1$
        \If{\( num_l > 10\)} ~~~\# \( ans_\mathbf{x} \) is not updated with ten continuous iterations.
            \State break
        \EndIf
    \EndWhile
    \State \Return Predicted solution \( ans_\mathbf{x} \); Equation residuals of $ans\_f$;
\end{algorithmic}
\end{algorithm}

\section{Numerical results}
In the following numerical experiments, we equip the two HANN methods with four hidden layers and 40 neurons per layer, and employ 1000 collocation points by the Latin hypercube sampling \cite{ms-1987} and the L-BFGS optimizer \cite{ln-1989} for the training, unless otherwise mentioned. Moreover, we assess the accuracy of predicted solutions by computing the system residual which are defined as the sum of absolute values of each equation residual. For brevity, we present an explicit procedure of HANN method for solving the first example and only show the final learning results for the other examples unless new techniques are introduced.  

\subsection{A single equation}
\begin{figure}[ht]
    \centering
    \begin{minipage}{0.33\textwidth}
    \centerline{\includegraphics[width=\linewidth]{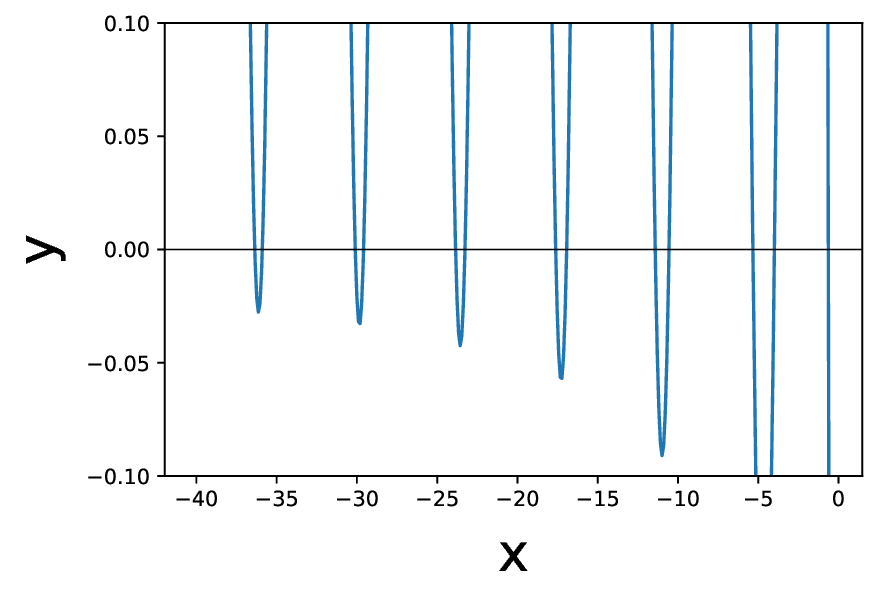}}
    \centerline{A}
    \end{minipage}
    \vspace{0.1cm}
    \begin{minipage}{0.3\textwidth}
    \centerline{\includegraphics[width=\linewidth]{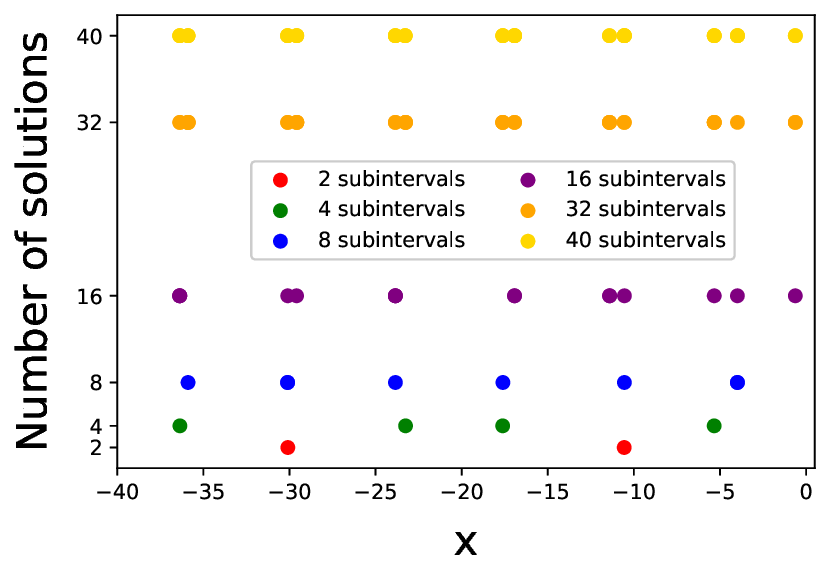}}
    \centerline{B}
    \end{minipage}
    \vspace{0.1cm}
    \begin{minipage}{0.33\textwidth}
    \centerline{\includegraphics[width=\linewidth]{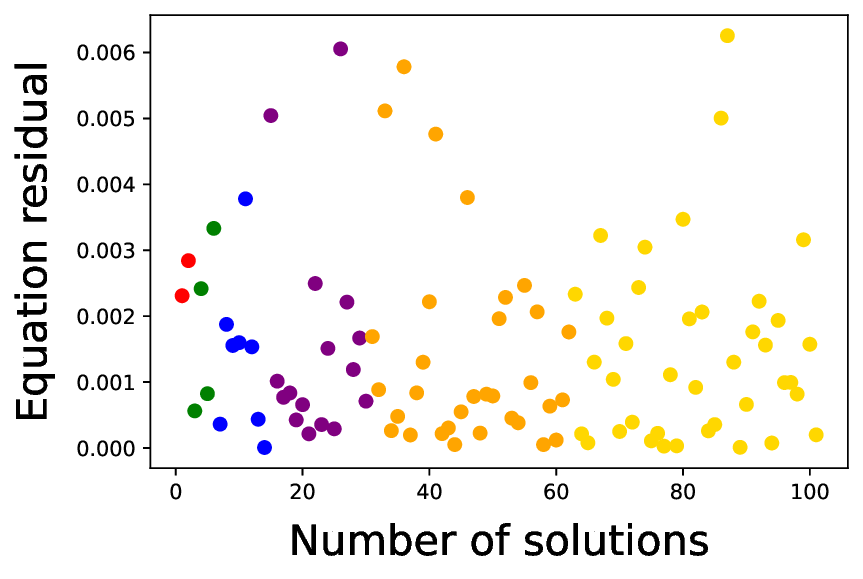}}
    \centerline{C}
    \end{minipage}
    \vspace{0.1cm}
    \begin{minipage}{0.4\textwidth}
    \centerline{\includegraphics[width=\linewidth]{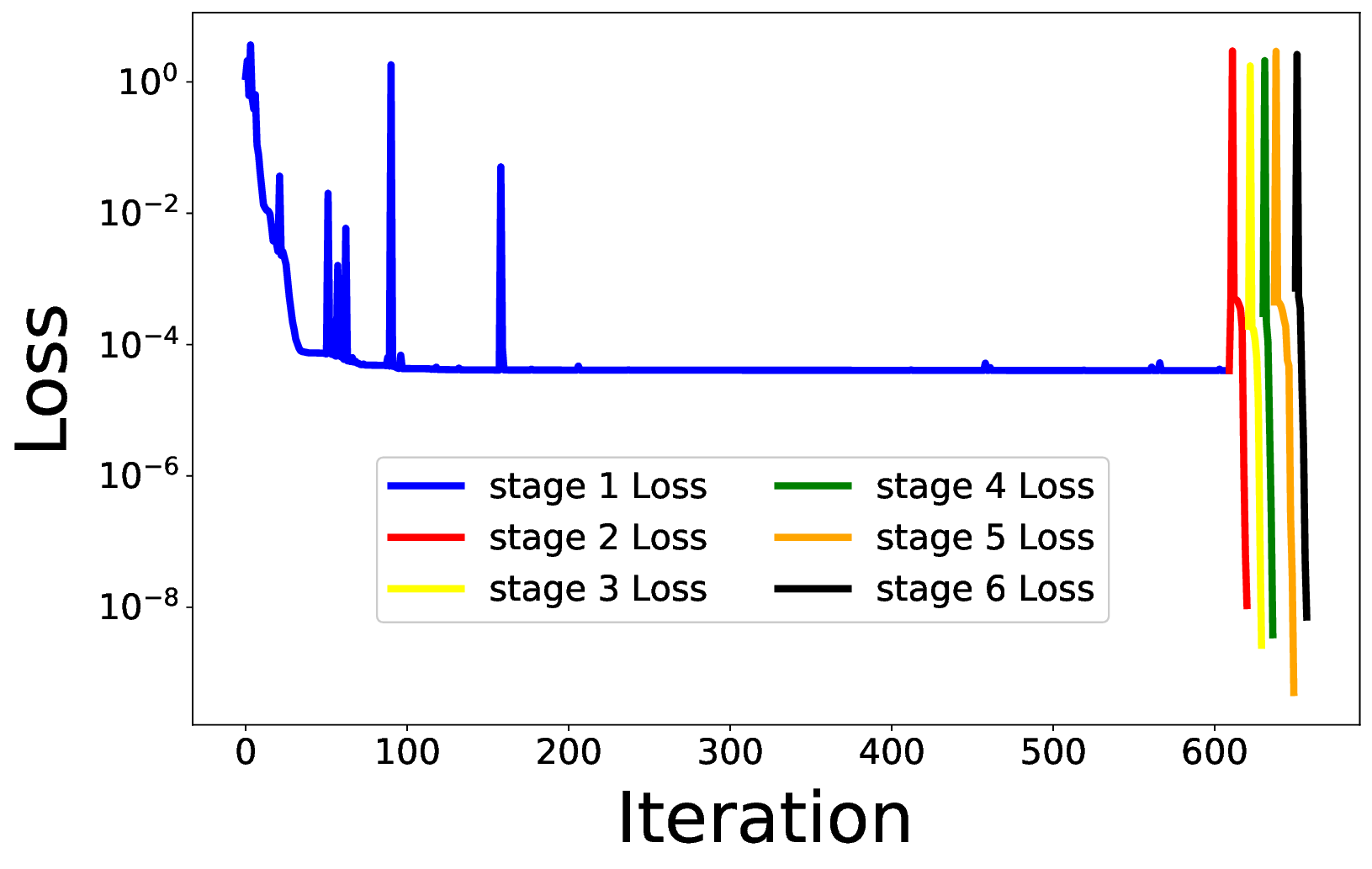}}
    \centerline{D}
    \end{minipage}
    \hspace{0.2cm} 
    \begin{minipage}{0.4\textwidth}
    \centerline{\includegraphics[width=\linewidth]{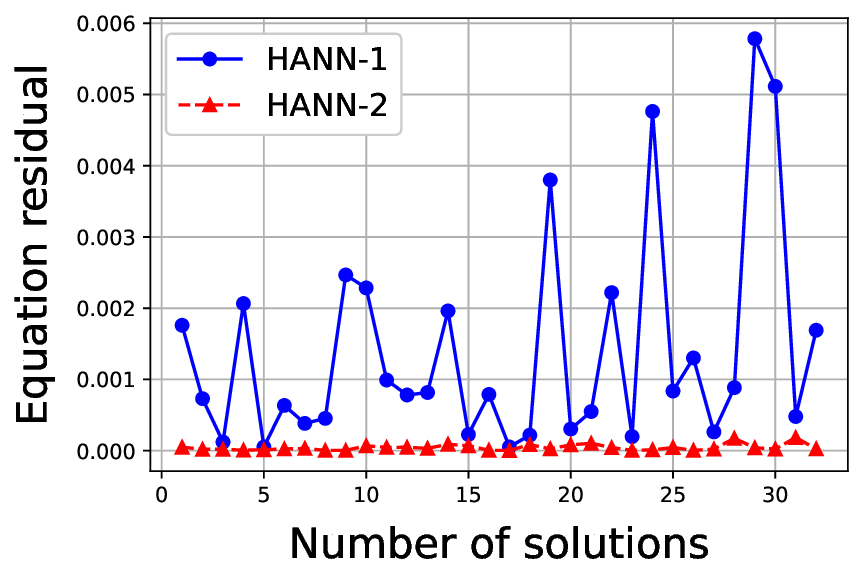}}
    \centerline{E}
    \end{minipage}
    \caption{\small (Color online) Equation (\ref{exa-1}): (A). Graph of $y = 1/x - \sin x + 1$. (B). Distribution of predicted solutions by HANN-1 under different numbers of sub-intervals. (C). Equation residual distribution of different predicted solutions in different sub-intervals, where the color of points correspond to the one of solutions in (B). (D). Loss history of HANN-2 with a random initial value -20.625. (E). Comparison of equation residuals of the learned 32 solutions by HANN-1 and HANN-2. }%
    \label{fig:sideex-1122}
\end{figure}

The first example is to consider a single equation \cite{4}
\begin{eqnarray}\label{exa-1}
&& f(x) = \frac{1}{x} - \sin x + 1 = 0,~~~ x\in (-40, 0),
\end{eqnarray}
which has 13 isolated solutions depicted by the graph of $y=f(x)$ in Figure \ref{fig:sideex-1122}(A). Next, we leverage the HANN-1 to find the 13 solutions by learning the solution $x(t)$ of the homotopy model
\begin{eqnarray}\label{exa-1-homo}
&& H(x,t) = t\left(\frac{1}{x} - \sin x + 1\right)+\gamma(1-t)\left(\frac{1}{x} - \sin x-\frac{1}{x_0} + \sin x_0\right) = 0,
\end{eqnarray}
where $x_0$ is an initial value, $t\in [0,1]$ denotes the homotopy parameter and $\gamma$ is the convergence-control parameter. Then, following the idea of homotopy continuation method \cite{li-1997}, the starting equation $H(x_0,0)=0$ holds identically and thus the initial condition $x(0)=x_0$ can be chosen arbitrary, and the final equation $H(x,1)=0$ corresponds to $f(x) =0$, thus after learning the solution $\hat{x}(\Theta; t)$ of model (\ref{exa-1-homo}), then $\hat{x}(\Theta; 1)$ is the required solution of equation (\ref{exa-1}), where $\Theta$ denotes the training parameter set. Then we employ the L-BFGS optimizer to minimize the loss function $\mathcal{L}(\Theta) = \mathcal{L}_{iv}(\Theta) + \mathcal{L}_h(\Theta)$ for model (\ref{exa-1-homo}), where
\begin{eqnarray} \label{eq:MSE-1}
\no &&\mathcal{L}_{iv}(\Theta)= \left|\hat{x}(\Theta; 0) - x_0 \right|^2, ~~~\mathcal{L}_h(\Theta)= \frac{1}{1000}\sum_{j=1}^{1000}\left|H(\hat{x}(\Theta; t_j),t_j) \right|^2,
\end{eqnarray}
and $\hat{x}(\Theta; t_j)$ denotes the output of HANN-1 and $t_j\in[0,1]$ is the input data. 

Before performing the training, we first determine an optimal value of the convergence-control parameter $\gamma$ via a systematic test by choosing the initial value $x_0=-15$, the seed $1234$ and 1000 collocation points. Consequently, Table \ref{tab:1-ex} shows that changing $\gamma$ from 5 to 0.01 enhances both training efficiency and solution accuracy, while decreasing $\gamma$ from 0.01 to 0.00001 yields negligible decreases on the computational costs but reduces the learned accuracy, thus if there is no special explanation, we fix $\gamma=0.01$ in here and the following experiments.
\begin{table}[htp]
\centering
\small
 \captionsetup{width=.65\textwidth}
\caption{Numerical results of different $\gamma$ on the convergence of HANN-1 for solving equation (\ref{exa-1}).}
\begin{tabular}{@{}cccccc@{}}
\toprule
$\gamma$  & $x$  & Equation residual&Training time(s) \\ \midrule
5	&-10.50394063	&2.323479e-02	&9.5779\\
1	&-17.66053462	&1.537179e-02	&10.9315\\
0.1	&-16.91877944	&4.990269e-03	&0.8091\\
0.01	&-17.61766674	&1.202379e-04	&0.6413\\
0.001	&-17.61837379	&3.578146e-04	&0.3714\\
0.0001	&-17.61923255	&6.470036e-04	&0.4382\\
0.00001	&-17.61945523	&7.221054e-04	&0.4117 \\\bottomrule
\end{tabular}
\label{tab:1-ex}
\end{table}

Next, in order to check the efficacy of HANN-1 with arbitrary initial value for solving equation (\ref{exa-1}), we first divide the defined interval (-40, 0) into 2, 4, 8, 16 and 32 equidistant sub-intervals in turn, and then choose midpoint in each sub-interval as the initial value to feed the HANN-1 for training. Figure \ref{fig:sideex-1122}(B) displays that, as the number of sub-intervals increases, the isolated solutions also appear consecutively, and with 32 sub-intervals, the HANN-1 gives the total 13 isolated solutions under the threshold $4.66\times10^{-2}$, where two learned solutions are regarded as one solution if the absolute value of their difference is smaller than the threshold. Moreover, even if the number of sub-intervals becomes bigger than 32, we still obtain the 13 solutions by HANN-1 under the same threshold, for example, the 40 sub-intervals in Figure \ref{fig:sideex-1122}(B).

In fact, Figure \ref{fig:sideex-1122}(C) exhibits that the HANN-1 can give an unique solution for each chosen initial value, but the distances among some solutions are very close and thus such solutions are regarded as one solution under the threshold. Therefore, the HANN-1 has the ability to find all isolated solutions of a single algebraic equation only if the number of equidistant sub-intervals is large enough. Moreover, the residuals of equation (\ref{exa-1}) inserted by the obtained solutions in Figure \ref{fig:sideex-1122}(B) are displayed in Figure \ref{fig:sideex-1122}(C) and stabilize at $10^{-3}$ order of magnitude, where the maximum is only $6.26\times10^{-3}$.
Thus, compared with the results by the Newton-Homotopy method \cite{4}, the HANN-1 can determine the total number of solutions within a fixed interval. 

Next, we employ the HANN-2 to further improve the accuracy of the learned solutions by HANN-1. Figure \ref{fig:sideex-1122}(D) displays the loss varying tendency of HANN-2 with a random initial value -20.62 over the iteration. In particular, stage one takes 609 iterations while in the other five stages the HANN-2 obviously accelerates the training procedure and only operates with total 48 iterations because the feeded initial values for the latter stage of HANN-2 are more closer to the true solution than the one of previous stage.
Furthermore, Figure \ref{fig:sideex-1122}(E) shows a comparison of equation residuals by HANN-1 and HANN-2 in the case of 32 sub-intervals where the red line for HANN-2 stabilizes in the neighbourhood of zero and is far below the blue line for HANN-1 which has big fluctuations, even in the worst case $1.83\times10^{-4}$ by HANN-2, there still has 61.77\% reduction of equation residual.
\subsection{Transcendental systems containing two equations}
In this section, we consider two systems of two transcendental systems which contain the absolute value function and the trigonometric function respectively. 
\subsubsection{System containing absolute value  function}
We consider a system of two equations containing the absolute value function and check the efficacy of HANN-1 for learning the solutions from a big defined domain \cite{gaa-2008}
\begin{eqnarray}\label{example-2}
\begin{cases}
 x^2 - y^2 = 0, \\
 1 - |x - y| = 0,~~~~~ (x, y)\in [-15, 15]\times[-15, 15],
 \end{cases}
\end{eqnarray}
whose two exact solutions are $(x,y)=(0.5,-0.5)$ and $(x,y)=(-0.5,0.5)$. Following the idea of HANN-1, we first construct a homotopy model for system (\ref{example-2}) defined by 
\begin{eqnarray}\label{example-2-homo}
&&\no  h_1(x,y,t):=t(x^2 - y^2)+\gamma(1-t)(x^2 - y^2-x_0^2 + y_0^2) = 0, \\
&&\no h_2(x,y,t):=t(1 - |x - y|) +\gamma(1-t)(|x_0 - y_0|- |x - y|)= 0,
\end{eqnarray}
where $t\in [0,1]$, $\gamma=0.01$ and $(x_0,y_0)$ is the initial value at $t=0$. Moreover, we divide the interval [-15, 15] of $x$ and $y$ into 7 equidistant grids respectively and use the 49 intersecting grid points as the initial value $(x_0,y_0)$, and then sample 1000 discretized points $t_j\in[0,1]$ as the input data for optimizing the loss function  $\mathcal{L}(\Theta) = \mathcal{L}_{iv}(\Theta) + \mathcal{L}_h(\Theta)$, where
\begin{eqnarray}
\no &&\mathcal{L}_{iv}(\Theta)= \left|\hat{x}(\Theta; 0) - x_0 \right|^2+\left|\hat{y}(\Theta; 0) - y_0 \right|^2, \\
&&\no \mathcal{L}_h(\Theta)= \frac{1}{1000}\sum_{j=1}^{1000}\left|h_1(\hat{x}(\Theta; t_j),\hat{y}(\Theta; t_j), t_j) \right|^2+\left|h_2(\hat{x}(\Theta; t_j),\hat{y}(\Theta; t_j), t_j) \right|^2,
\end{eqnarray}
and $\hat{x}(\Theta; t_j)$ and $\hat{y}(\Theta; t_j)$ are the two outputs of HANN-1.

Since system (\ref{example-2}) has two distinct solutions, we first check the efficacy of HANN-1 with different initial values and different random seeds. Table \ref{tab:results1} shows that the HANN-1 with either different random seeds and the same initial value  or the same random seed and different initial values generates different solutions. 
Furthermore, for the initial values $(0,0)$ and $(-5,-5)$ which locate at the equidistant position of the two true solutions, the HANN-1 with the two random seeds generate two different solutions of system (\ref{example-2}). Whereas, for the initial values $(5,5)$ equidistant with the two true solutions, the HANN-1 with two different random seeds gives different solutions. More interesting, for the initial value $(-5,5)$ which is closer to the reference solution $(-0.5,0.5)$, but the HANN-1 predicts the solution approximating $(0.5,-0.5)$. Similarly, the HANN-1 with initial value $(5,-5)$ also gives the solution near $(0.5,-0.5)$ but not $(-0.5,0.5)$. Such prediction results mean that the HANN-1 is not sensitive to the position of initial value and different with the Newton method which is strongly dependent on the selection of initial value.
 \begin{table}[htp]
\centering
\small
 \captionsetup{width=.75\textwidth}
\caption{\small Numerical results by HANN-1 for system (\ref{example-2}) with different initial values and different random seeds. }
\label{tab:results1}
\begin{tabular}{@{}cccccc@{}}
\toprule
random seed and initial value & $x$ & $y$ & Equation residual \\ \midrule
\hspace{0.3cm}1\hspace{0.3cm} (0, 0) & 0.50168678&	-0.50079411&	3.375783e-03\\
1234 (0, 0)&-0.49992005	&0.50106249	&2.126111e-03\vspace{4pt}\\
\hspace{0.3cm}1\hspace{0.3cm} (-5, 5)&	-0.49983449&	0.49940003&	1.199607e-03\\
1234	(-5, 5)	&0.50019711	&-0.49749038	&5.012977e-03\vspace{4pt}\\
\hspace{0.3cm}1\hspace{0.3cm} (5, -5)	&0.49686094	&-0.49906324	&6.269144e-03\\
1234 (5, -5)	&-0.49860885	&0.49963823&2.780499e-03\vspace{4pt}\\
\hspace{0.3cm}1\hspace{0.3cm} (5, 5)	&-0.49957557	&0.49940736&	1.185112e-03\\
1234	(5, 5)	&0.50001307&	-0.50172158&	3.446114e-03\vspace{4pt}\\
\hspace{0.3cm}1\hspace{0.3cm}(-5, -5)&	-0.50065917&	0.50010156&	1.318759e-03\\
1234 (-5, -5)&	-0.50142112&	0.50019246&	2.844226e-03\\
\bottomrule
\end{tabular}
\end{table}

\begin{figure}[ht]
    \centering
        \begin{minipage}{0.33\textwidth}
    \vspace{-0.5cm}
    \centerline{\includegraphics[width=\linewidth]{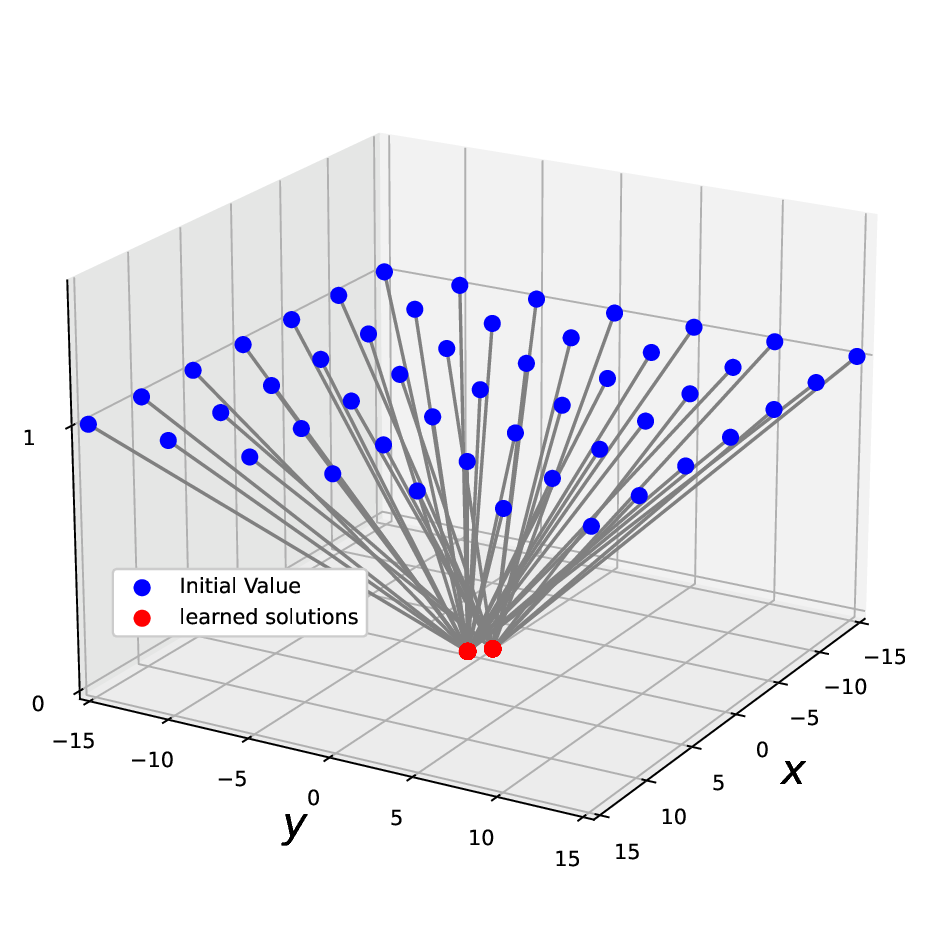}}
    \centerline{A}
    \vspace{0.1cm}
    \end{minipage}
     \hspace{0.5cm}
    \begin{minipage}{0.45\textwidth}
    \centerline{\includegraphics[width=\linewidth]{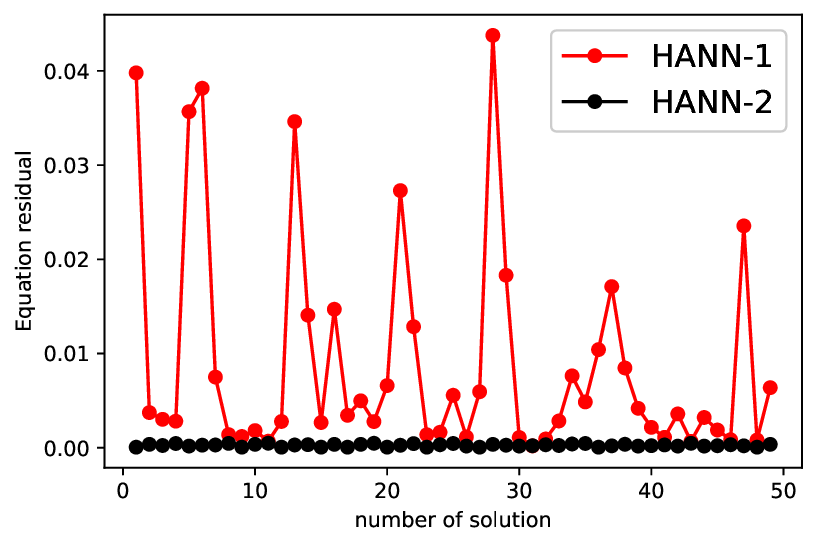}}
    \centerline{B}
    \vspace{0.1cm}
    \end{minipage}
\begin{minipage}{0.32\linewidth}
    \centerline{\includegraphics[width=1.0\textwidth]{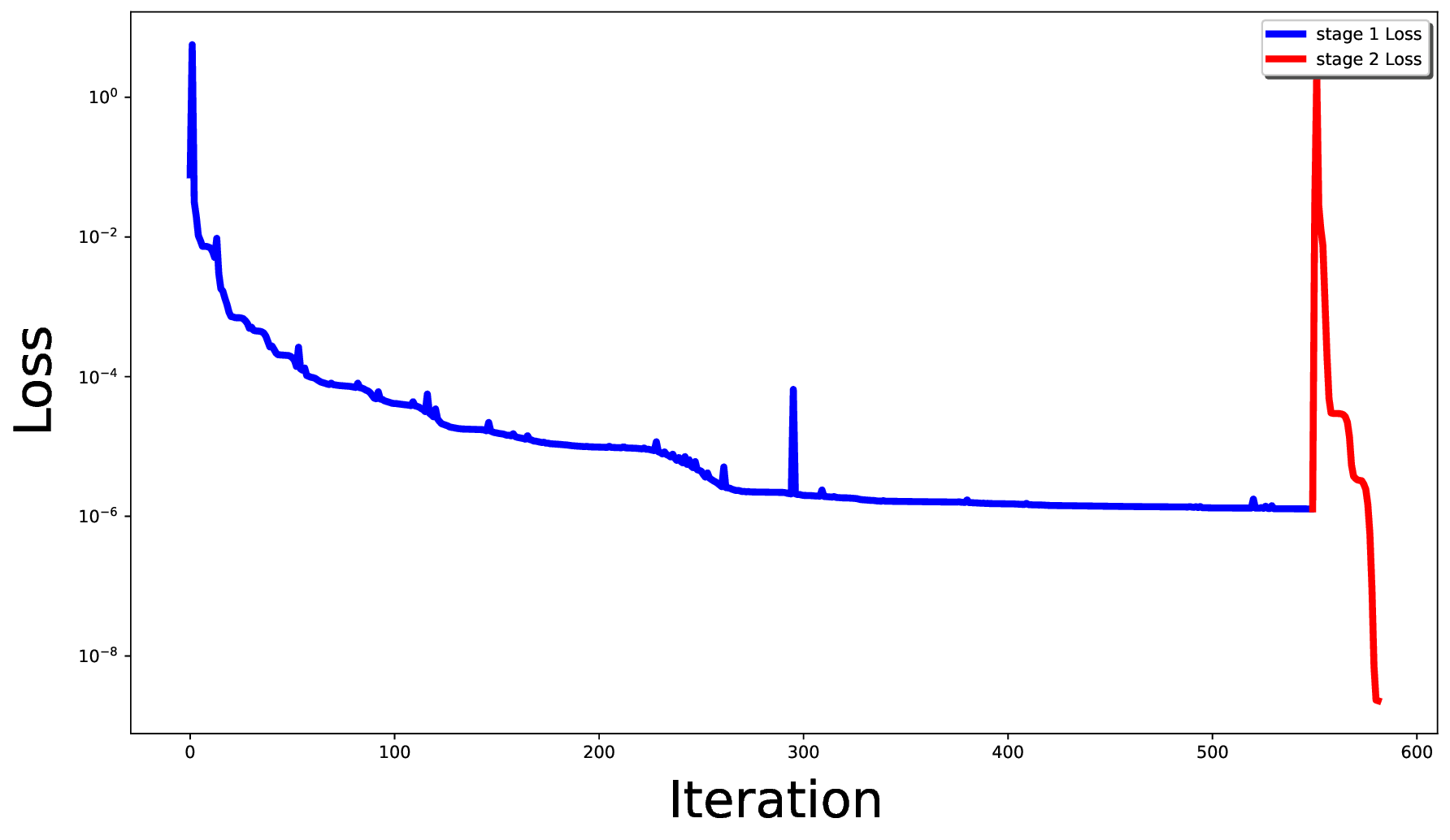}}
    \centerline{C}
	\end{minipage}
	\begin{minipage}{0.3\linewidth}
    \centerline{\includegraphics[width=1.0\textwidth]{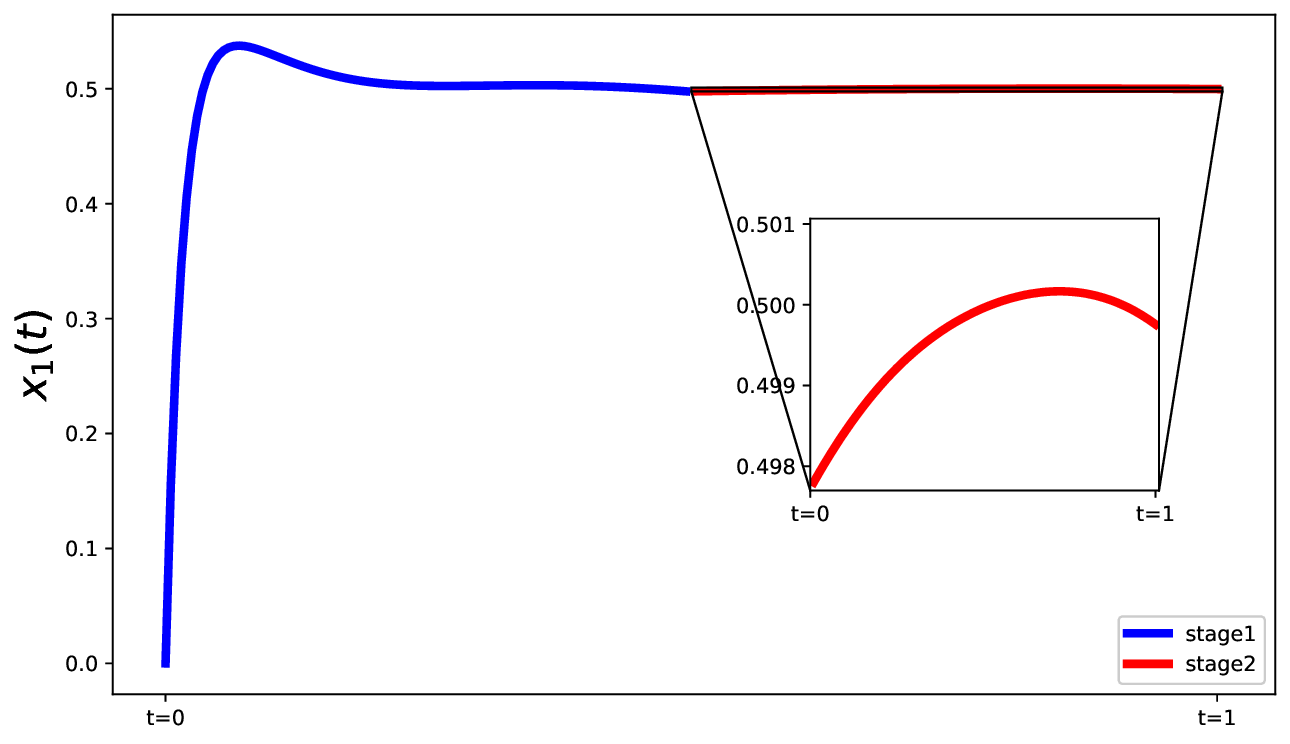}}
    \centerline{D}
	\end{minipage}
	\begin{minipage}{0.3\linewidth}
    \centerline{\includegraphics[width=1.0\textwidth]{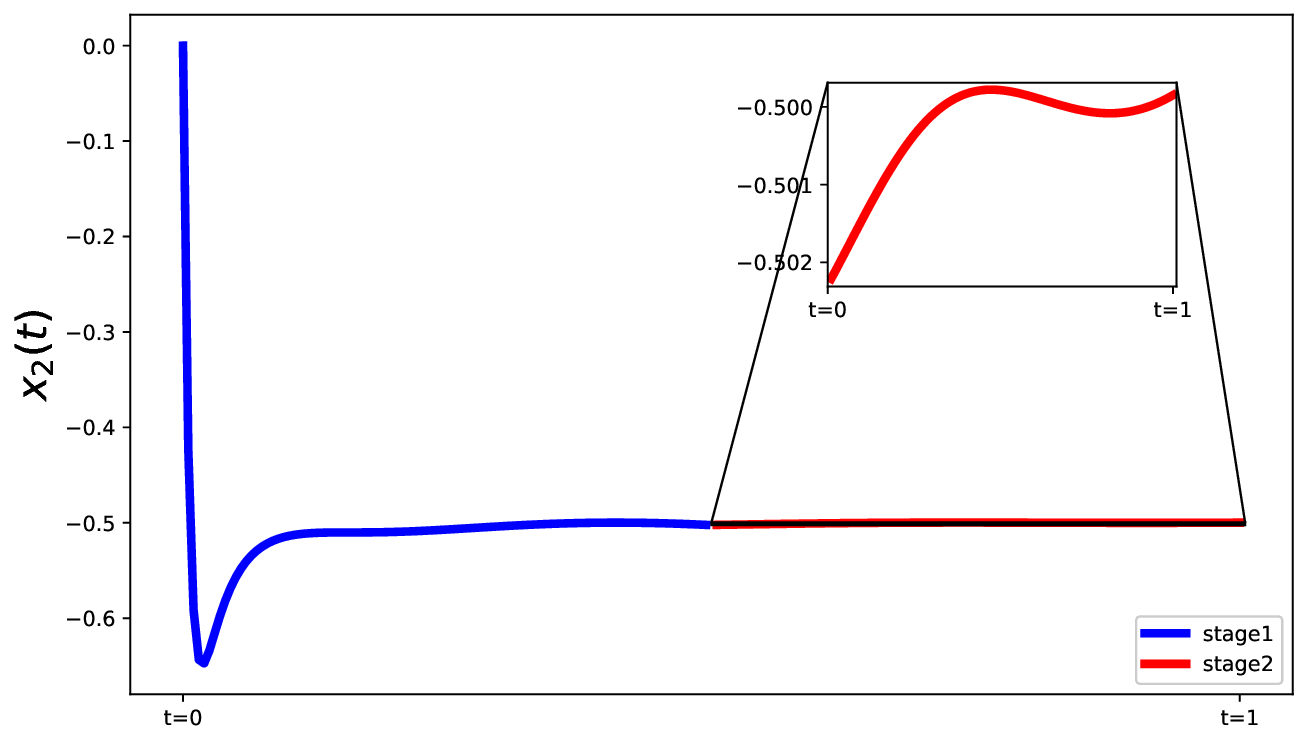}}
    \centerline{E}
	\end{minipage}
    \caption{\small(Color online) System (\ref{example-2}): (A). The mappings by HANN-1 from the 49 initial values to two predicted solutions where the two solid red points are the predicted solutions under the threshold \(3.54\times 10^{-2}\); (B). Comparisons of equation residuals of the 49 predicted solutions with HANN-1 and HANN-2. (C) Loss decent tendency by HANN-2 over iterations. (D) Homotopy curve $x_1(t)$. (E) Homotopy curve $x_2(t)$.  }%
    \label{fig:sideex-222a}
\end{figure}

Furthermore, Figure \ref{fig:sideex-222a}(A) shows that HANN-1 with 49 initial values gives the two expected solutions of system  (\ref{example-2}) visible to the naked eye and their accuracies are further improved by HANN-2 in Figure \ref{fig:sideex-222a}(B),  where the equation residuals for the 10 of 49 solutions with HANN-2 arrive at $10^{-5}$ order of magnitude, while the maximal residual of other 39 solutions is \(4.82\times 10^{-4}\) which still has 30.95\% improvement than the corresponding one by HANN-1. Moreover, Figure \ref{fig:sideex-222a}(C) shows the loss history of HANN-2 with the initial value $(0, 0)$ where stage two is obviously faster than stage one and the minimum arrives at $2.19\times10^{-9}$, while Figure \ref{fig:sideex-222a}(D, E) shows the homotopy curves of $x_1(t)$ and $x_2(t)$ over the iterations by HANN-2, which are continuous across the whole interval.  Consequently, the learned solutions in stage one which locate at the right ends of solid blue curve, $x_1(1)=0.497794$ and $x_2(1) = -0.502212$, and their relative errors to the true solution (0.5, -0.5) are $Re(x_1) = 0.4412\%, Re(x_2) = 0.4424\%$. By contrast, the predicted solutions of stage two locating at the right ends of red curve, $x_1(1)=0.499757$ and $x_2(1) = -0.499849$, and the corresponding relative errors are $Re(x_1) = 0.0486\%, Re(x_2) = 0.0302\%$, which has one order of magnitude improvement than stage one.

\subsubsection{System containing trigonometric functions}
The other transcendental system contains the trigonometric function  and is used to evaluate the capability of HANN-1 to search for isolated solutions in the finite domain
\begin{equation}\label{example-3}
\begin{cases}
2(x_2 - x_1) + \sin(2x_2) - \sin(2x_1) - 1.2 = 0,\\
\cos(2x_1) - \cos(2x_2) - 0.4 = 0, ~~~(x_1,x_2)\in [-5,5]\times[-5,5],
\end{cases}
\end{equation}
where $x_1$ and $x_2$ describe the joint angles between mechanical arms in inverse kinematics \cite{18}.

\begin{figure}[ht]
    \centering
    \begin{minipage}{0.33\textwidth}
    \vspace{-0.1cm}
    \centerline{\includegraphics[width=\linewidth]{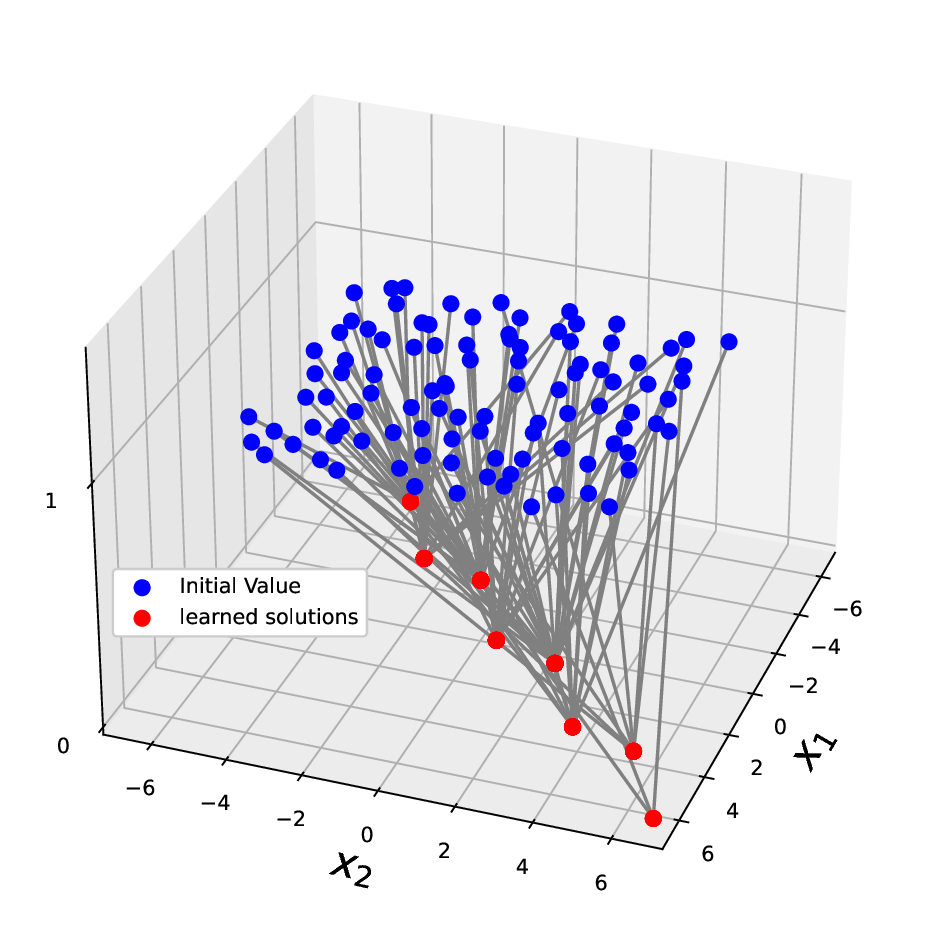}}
    \centerline{A}
    \end{minipage}
     \hspace{0.5cm}
    \begin{minipage}{0.45\textwidth}
    \centerline{\includegraphics[width=\linewidth]{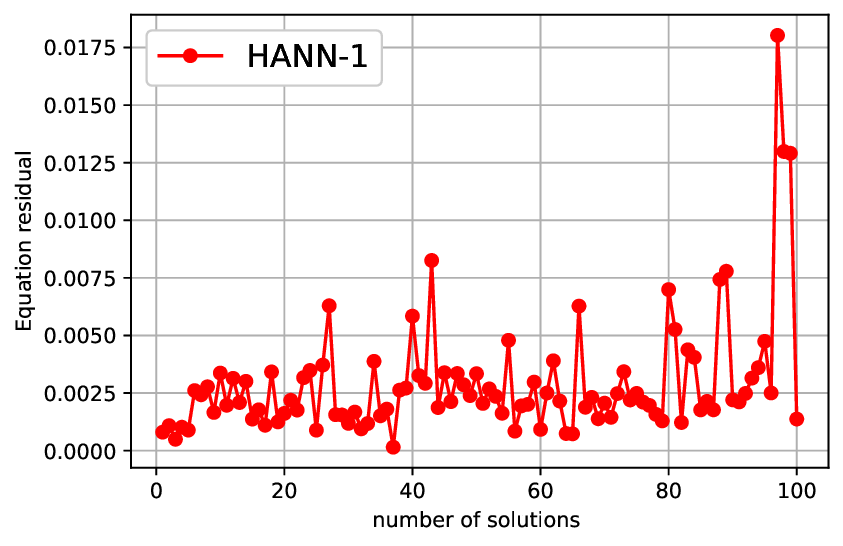}}
    \centerline{B}
    \end{minipage}
    \caption{\small(Color online) System (\ref{example-3}): (A). The mappings by HANN-1 with the 100 initial values where the eight accumulated solid red points are the predicted solutions; (B). Equation residuals of the predicted solutions by HANN-1 with the 100 initial values.}%
    \label{fig:sideex-222}
\end{figure}

In order to systematically evaluate the efficacy of HANN-1 with different initial values, we first discretize the interval $[-5,5]$ for $x_1$ and $x_2$ into ten equidistant grids respectively and then obtain one hundred squares with area one on the plane, and finally randomly choose one point in each square as the initial value of HANN-1. 
Then Figure \ref{fig:sideex-222}(A) shows that HANN-1 with the 100 random initial values gives eight isolated solutions of system (\ref{example-3}) described by the solid red points under the threshold $1.08\times10^{-2}$, while Figure \ref{fig:sideex-222}(B) displays the equation residuals of the obtained 100 solutions where, except for the second to the fourth solutions from the bottom, the others all arrive at $10^{-3}$ order of magnitude.

In particular, Table  \ref{tab:example3.2} displays the specific eight solutions of system (\ref{example-3}) by HANN-1 and shows a comparison with the results by the five known methods \cite{18}. Obviously, the HANN-1 not only recovers the known solution but also gives seven new high-accuracy solutions where four solutions fall in the defined domain and three solutions are beyond the domain, while the HANN-2 further refines the learned solutions by HANN-1 and presents one order of magnitude improvement in most cases. Note that the residuals of two equations in system (\ref{example-3}) are listed separately in order to facilitate the comparisons with the results of five known methods.
\begin{table}[ht]
 \captionsetup{width=.75\textwidth}
\centering
\small
\caption{\small Numerical solutions of system (\ref{example-3}) with seven methods, where the eight solutions by HANN-1 are obtained with 100 random initial values  under the residual threshold $1.08\times10^{-2}$ while the initial values of HANN-2 are the learned eight solutions by HANN-1.}
\label{tab:example3.2}
\renewcommand{\arraystretch}{1.1}
\begin{tabular}{@{}ccccccc@{}}
\toprule
Method & Solution $(x_1, x_2)$ & Equation residual \\ \midrule
Newton's method&(0.15, 0.49)&(-1.68e-03, 1.5e-02)\\
Secant method&(0.15, 0.49)&(-1.68e-03, 1.5e-02)\\
Broyden’s method&(0.15, 0.49)&(-1.68e-03, 1.5e-02)\\
Effati’s method&(0.1575, 0.4970)&(5.46-03, 7.39e-03)\\
Evolutionary approach&(0.15772, 0.49458)&(1.26e-03, 9.69e-04)\\ \hline
&(0.15404579, 0.49054697)	&(3.20e-03, 8.68e-04) \\
&(-2.9850282, -2.648330057)&(2.17e-04, 5.26e-04)\\
&(-2.46140343, -0.881259401)&(5.81e-04, 6.67e-04)\\
&(0.680151726, 2.259743585)&(6.50e-04, 1.96e-04)\\
HANN-1&(3.29719883,  3.63507387)&(7.36e-04, 3.89e-04)\\
&(3.822020409, 5.401734489)&(6.72e-04, 2.20e-04)\\
&(-5.603113604, -4.023228239)&(3.87e-04, 3.42e-04)\\
&(6.440921995, 6.77761628)&(1.01e-03, 1.48e-03)\\
\hline
& (0.15680684, 0.49370563)& (3.73e-04, 9.77e-05)  \\
&(-2.98500176, -2.64813161)&(9.76e-05, 1.37e-05)\\
&(-2.46132136, -0.88150516)&(2.59e-04, 7.07e-05)\\
&(0.68026611, 2.26008618)&(2.47e-04, 8.11e-05)\\
HANN-2&(3.29856214, 3.63538272)&(4.13e-04, 4.71e-04)\\
&(3.82186076, 5.40167894)&(2.51e-04, 7.65e-05)\\
&(-5.60291107, -4.0230971)&(2.67e-04, 6.48e-05)\\
&(6.44014078, 6.77695557)&(3.89e-04, 4.77e-04)\\
\bottomrule
\end{tabular}
\end{table}

In addition, we employ the Python's built-in Fsolve function with either the initial values in Table \ref{tab:example3.2} obtained by HANN-1 or the ones directly feeding into HANN-1 for solving systems (\ref{example-2}) and (\ref{example-3}) respectively, and find that the Fsolve function fails to work. Therefore, the HANN-1 has strong capability to solve the transcendental equations.
\subsection{High-dimensional systems containing ten equations}

Observed that the HANN-1 presents excellent jobs for the single equation and the two transcendental systems of two equations, we next extend it to solve the high-dimensional systems containing ten equations which include an ill-condition system.
\subsubsection{Interval arithmetic benchmark}
\begin{figure}[ht]
    \centering
    \begin{minipage}{0.45\textwidth}
    \centerline{\includegraphics[width=\linewidth]{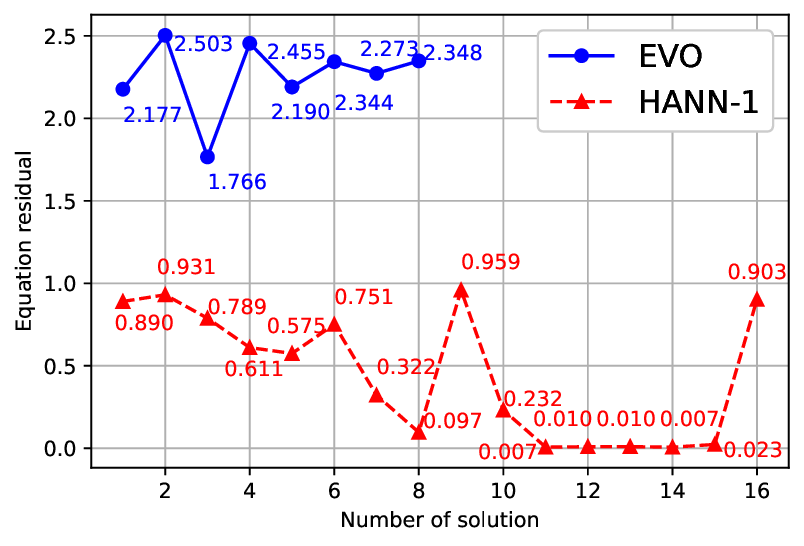}}
    \centerline{A}
    \end{minipage}
    \hspace{0.2cm} 
    \begin{minipage}{0.45\textwidth}
    \centerline{\includegraphics[width=\linewidth]{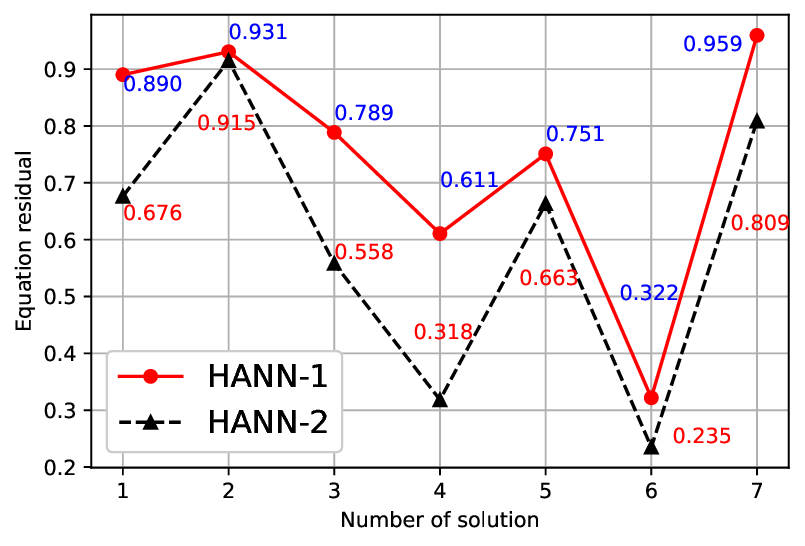}}
    \centerline{B}
    \end{minipage}
    \caption{\small(Color online) System (\ref{example-6}): Comparisons of equation residuals among the EVO and HANN-1 and HANN-2 methods: (A) EVO and HANN-1 (B) HANN-1 and HANN-2 for the seven solutions than cannot be improved the Fsolve function.}%
    \label{fig:sideex-1}
\end{figure}

We first check the efficacy of HANN-1 for solving the normal multiple-variables system by considering a benchmark problem from interval arithmetic \cite{29}
\begin{equation}\label{example-6}
\begin{cases}
x_1 - 0.25428722 - 0.18324757\, x_4x_3x_9 =0,\\
x_2 - 0.37842197 - 0.16275449\, x_1x_{10}x_6 =0, \\
x_3 - 0.27162577 - 0.16955071\, x_1x_2x_{10}  =0,\\
x_4 - 0.19807914 - 0.15585316\, x_7x_1x_6  =0,\\
x_5 - 0.44166728 - 0.19950920\, x_7x_6x_3  =0,\\
x_6 - 0.14654113 - 0.18922793\, x_8x_5x_{10} =0, \\
x_7 - 0.42937161 - 0.21180486\, x_2x_5x_8  =0,\\
x_8 - 0.07056438 - 0.17081208\, x_1x_7x_6  =0,\\
x_9 - 0.34504906 - 0.19612740\, x_{10}x_6x_8  =0,\\
x_{10} - 0.42651102 - 0.21466544\, x_4x_8x_1 =0,
\end{cases}
\end{equation}
which is a system of three degree polynomials with ten variables and considered as illustrated examples in many literatures \cite{27,2}.  Next, we first artificially define the domain of $(x_1,x_2,\dots,x_{10})\in \Omega=\prod_{i=1}^{10} [-30,30]$ and employ the Latin hypercube sampling to get 200 grid points in $\Omega$, then choose them as the initial values of HANN-1 instead of training points for learning solutions of system (\ref{example-6}). In particular, for the sake of exhibiting strong learning capability of HANN-1 for dealing with normal system (\ref{example-6}), we equip HANN-1 with an extremely simple architecture, two hidden layers and two neurons in each layer, and leverage five collocation points and convergent-control parameter $\gamma=0.0001$ for the training. Consequently, we obtain sixteen solutions with accuracy arriving at $10^{-1}$ order of magnitude and show them in Figure \ref{fig:sideex-1}(A), where the dotted red line describing the equation residuals of sixteen solutions by HANN-1 stays far below the solid blue line by the evolutionary algorithm (EVO) in \cite{2}, at least one order of magnitude improvement and three orders in the best case.

In particular, we employ the built-in Fsolve function in Python to further improve the accuracies of sixteen solutions, where only nine solutions are refined and shown in Table \ref{tab:results}, the other seven solutions that cannot be refined by the Fslove function are successfully enhanced by HANN-2 and are displayed in Figure \ref{fig:sideex-1}(B). More interesting, by directly applying the aforementioned 200 random grid points as the initial values of the built-in Fsolve function in Python, instead of using the HANN-1 initialization, we only obtain seven solutions with accuracy higher than $10^{-1}$ where the five high-accuracy solutions present the equation residuals $1.05\times 10^{-9}$, $3.28\times 10^{-9}$, $2.89\times 10^{-13}, 6.97\times 10^{-9}$ and $4.64\times 10^{-10}$, whose accuracies are lower than their counterparts of the $k$-th ($k=1,2,6,7,8$) case in Table \ref{tab:results}. Therefore, the HANN-1 can rapidly give a candidate set of rough learned solutions, then collecting HANN-1 with the Fsolve function or HANN-2 may generate new high-accuracy solutions.
\begin{table}[ht]
\centering
\footnotesize
 \captionsetup{width=.85\textwidth}
\caption{\small Numerical results of system (\ref{example-6}) by the Fslove function where the initial values are the learned solutions by HANN-1.}
\begin{tabular}{@{}cccccc@{}}
\toprule
Case &$x_1$ & $x_2$ & $x_3$ & $x_4$ & $x_5$  \\\midrule
1&-2.412220977&	-2.290323698&	-2.111011823&	-2.221227578&	-2.269731897\\
2&-2.068464124&	2.358431355&	2.102111246&	2.397098585&	-2.418618802\\
3&-0.017212478&	0.410907113&	0.367100602&	10.06725143&	-269.1777816\\
4&-0.014075949&	0.397708215&	0.298453759&	11.78783966&	-314.0994619\\
5&-0.006496133&	0.358320826&	0.157988754&	14.77312841&	-453.394766\\
6&0.257431972&	0.380809764&	0.279386975&	0.200718638&	0.444888448\\
7&1.843705186&	1.969576487&	1.620416551&	2.085091918&	2.562787287\\
8&2.06179032&	-1.864657949&	-1.402607323&	-2.033457959&	2.385743371\\
9&2.420875356&	-1.961971859&	-1.937699599&	-1.998564649&	2.3727384\\

\bottomrule
 $x_6$ & $x_7$ & $x_8$ & $x_9$ & $x_{10}$& Equation residual \\\midrule
-2.66909418&	-2.412273916&	-2.581165424&	-3.098083885&	-2.544602352&	9.36E-11\\
2.657862106&	-2.566873252&	2.48082643&	-2.515310589&	-2.212241962&	2.39E-13\\
14.44662791&	-254.8257341&	10.89577516&	-0.430031337&	-0.025758793&	3.88E-08\\
15.64990574&	-337.5436885&	12.77344875&	-0.444130318&	-0.020424264&	6.14E-07\\
26.08755186&	-551.9110857&	16.05172001&	-1.109535738&	-0.018833827&	1.31E-08\\
0.145833685&	0.430201613&	0.073404607&	0.346064111&	0.427182476	&1.41E-14\\
2.419666215&	2.716139943&	2.138655983&	2.569032396&	2.191674122	&2.74E-12\\
-2.604983832&	2.666905401&	-2.375158355&	3.458373215&	2.56712054	&6.64E-13\\
-2.276614347&	2.395151577&	-2.105137056&	2.903343349&	2.643858443	&2.96E-09\\
\bottomrule
\end{tabular}
\label{tab:results}
\end{table}

\subsubsection{Combustion problem}
 We consider the combustion chemistry example at a temperature of $3000^{\circ}$C proposed in \cite{35,50}, which is described by the following sparse system of equations
\begin{equation}\label{example-7}
\begin{cases}
x_2 + 2x_6 + x_9 + 2x_{10} = 10^{-5}, \\
x_3 + x_8 = 3 \cdot 10^{-5}, \\
x_1 + x_3 + 2x_5 + 2x_8 + x_9 + x_{10} = 5 \cdot 10^{-5}, \\
x_4 + 2x_7 = 10^{-5}, \\
0.5140437 \cdot 10^{-7} x_5 = x_1^2, \\
0.1006932 \cdot 10^{-6} x_6 = 2x_2^2, \\
0.7816278 \cdot 10^{-15} x_7 = x_4^2, \\
0.1496236 \cdot 10^{-6} x_8 = x_1 x_3, \\
0.6194411 \cdot 10^{-7} x_9 = x_1 x_2, \\
0.2089296 \cdot 10^{-14} x_{10} = x_1 x_2^2.
\end{cases}
\end{equation}

System (\ref{example-7}) is highly possible an ill-conditioned equation since there exist huge differences among the coefficients of system (\ref{example-7}), even arriving at $10^{15}$ times in the last equation, which may lead to poor accuracies and big deviations with true solutions. Comparing with the normal system (\ref{example-6}), systematical investigations under different initial values and random seeds in Table \ref{tab:4.4} verify that the HANN-1 for system (\ref{example-7}) is sensitive to both the initial value and the random seed. 
In particular, for the first four cases of same initial value and different seeds, the HANN-1 gives completely different solutions. Meanwhile, in the latter four cases of different initial values and same seeds the HANN-1 still generates different solutions whose corresponding equation residuals arrive at $10^{-2}$ order of magnitude, significantly better than the evolutionary algorithm \cite{2}. Thus, the numerical results validate that the HANN-1 has the capability to deal with the ill-conditioned nonlinear algebraic system.
\begin{table}[ht]
 \captionsetup{width=0.95\textwidth}
\centering
\small
\footnotesize
\caption{\small Numerical results of system (\ref{example-7}) by HANN-1 under different random seeds and different initial values.}
\begin{tabular}{@{}ccccccccccc@{}}
\toprule
 Random seed and Initial Value&$x_1$ & $x_2$ & $x_3$ & $x_4$ & $x_5$  \\\midrule
1	\hspace{0.4cm}[0,0,0,0,0,0,0,0,0,0]&	-0.00370259&	0.01767804&	-0.13773655&	-0.04282948&	-0.10548795\\
12	\hspace{0.26cm}[0,0,0,0,0,0,0,0,0,0]&	0.03296258&	-0.03012133&	0.03902607&	-0.01410685&	0.10905272\\
123\hspace{0.13cm}	[0,0,0,0,0,0,0,0,0,0]&	-0.00646666&	-0.00452143&	-0.06526272&	0.01582104&	-0.01067824\\
9999	[0,0,0,0,0,0,0,0,0,0]&	-0.00878543&	0.01154406&	-0.07431326&	-0.02868609&	0.00917935 \vspace{4pt} \\
1234	[0,0,0,0,0,0,0,0,0,0]&	0.00716163&	-0.00465711&	0.13905826&	0.03587913&	-0.00408339\\
1234	[1,1,1,1,1,1,1,1,1,1]&	-0.00247197&	-0.02307291&	0.90155323&	-0.01031309&	0.69216196\\
1234	[2,2,2,2,2,2,2,2,2,2]&	-0.00017406&	-0.0000842&	2.01570322&	-0.02618407&	1.47928839\\
1234	[-1,-1,-1,-1,-1,-1,-1,-1,-1,-1]&	0.00805636&	-0.03060408&	0.95590331&	-0.03854047&	-0.67681669\\
\bottomrule
 $x_6$ & $x_7$ & $x_8$ & $x_9$ & $x_{10}$& Equation residual \\\midrule
-0.09386463&	0.02095912&	0.13719451&	-0.01115377&	0.08632501&	1.682803e-02\\
0.11717426&	0.00811038&	-0.04010447&	-0.22151589&	0.00964912&	1.327700e-02\\
0.03911695&	-0.00977146&	0.06610348&	-0.01062654&	-0.03195714&	9.948402e-03\\
0.14285328&	0.01244537&	0.0769432&	0.12151023&	-0.20886281&	1.046578e-02\vspace{4pt}\\
-0.02581832&	-0.02078766&	-0.14286198&	0.23417293&	-0.08840915&	1.576982e-02\\
0.98678727&	0.00684249&	-0.90100766&	0.9914312&	-1.471531&	8.931228e-03\\
1.94125993&	0.01338621&	-2.01682465	&2.00679399&	-2.94418104&	6.449540e-03\\
-0.88533855&	0.01934219&	-0.95308232&	2.80243835&	-0.50116689&	2.010216e-02\\
 \bottomrule
\end{tabular}\label{tab:4.4}
\end{table}

\begin{figure}[htp]
    \centering
    \begin{minipage}{0.45\linewidth}
		\vspace{3pt}
		\centerline{\includegraphics[width=\textwidth]{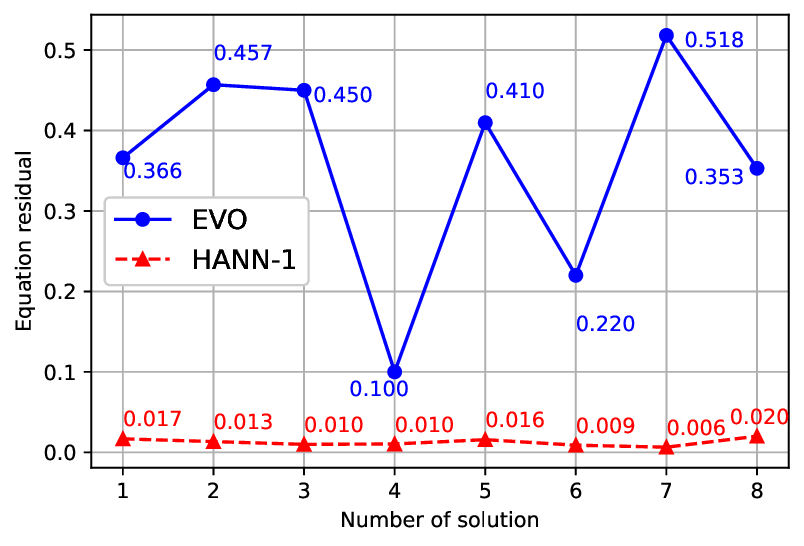}}
        \centerline{A}
	\end{minipage}
\hspace{0.2cm}
    \begin{minipage}{0.45\linewidth}
		\vspace{3pt}
		\centerline{\includegraphics[width=\textwidth]{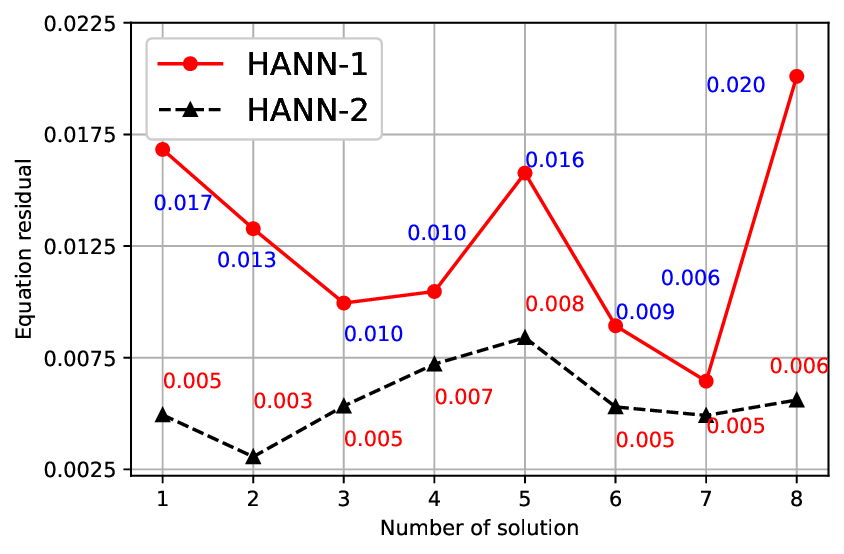}}
        \centerline{B}
	\end{minipage}
    \caption{\small(Color online) System (\ref{example-7}): Comparisons among the EVO, HANN-1 and HANN-2 methods. (A) EVO and HANN-1 (B) HANN-1 and HANN-2.}
    \label{fig:25}
\end{figure}

Furthermore, Figure \ref{fig:25} shows the performances of three algorithms in solving system (\ref{example-7}) respectively. It is evident that HANN-1 has better learning results than the EVO \cite{2} in Figure \ref{fig:25}(A) where the worst case by HANN-1 still has one order of magnitude improvement of equation residual. By contrast, the HANN-2 exhibits better performances than the HANN-1 in Figure \ref{fig:25}(B), where the worst case occurs for the seventh solution which still has $17\%$ improvement of equation residual. In particular, eight solutions are found without preprocessing the system, such as the usual scaling transformations $x_i=10^{-5}z_i\,(i=1,2,\dots,10)$. 
\subsection{Multivariate time-varying equations}
The key point of HANN is to artificially introduce a homotopy model of the algebraic equations to be solved in which the homotopy parameter $t\in[0,1]$ is chosen as the input variable of neural network and make the training work. It is worthy of noting that the time-varying problems themselves possess the time variable and thus naturally are suitable for the framework of HANN. Therefore, we extend the idea of HANN-1 to study the time-varying problems by considering a representative time-varying nonlinear system
\begin{equation}\label{example-8}
\begin{cases}
f_1:= \ln x_1(t) - \frac{1}{t+1} = 0, \\
f_2:= x_1(t)x_2(t) - \exp\big(\frac{1}{t+1}\big)\sin t = 0, \\
f_3:= x_1^2(t) - \sin(t)x_2(t) + x_3(t) - 2 = 0, \\
f_4:= x_1^2(t) - x_2^2(t) + x_3(t) + x_4^2(t) - t = 0,~~~t\in[0,10],
\end{cases}
\end{equation}
which has been used to check the efficacy of many proposed methods \cite{30,xzl-2019}. Moreover, for the verification purposes, we present the exact solution by
\begin{eqnarray}\label{true-so}
 x_1^*(t)=\exp\big(\frac{1}{t+1}\big), ~ x_2^*(t)=\sin t, ~ x_3^*(t)=2 - \exp\big(\frac{2}{t+1}\big) + \sin^2 t,~ x_4^*(t)=t-2.
\end{eqnarray}

\begin{figure}[htp]
    \centering
    \begin{minipage}{0.46\linewidth}
		\vspace{3pt}
		\centerline{\includegraphics[width=\textwidth]{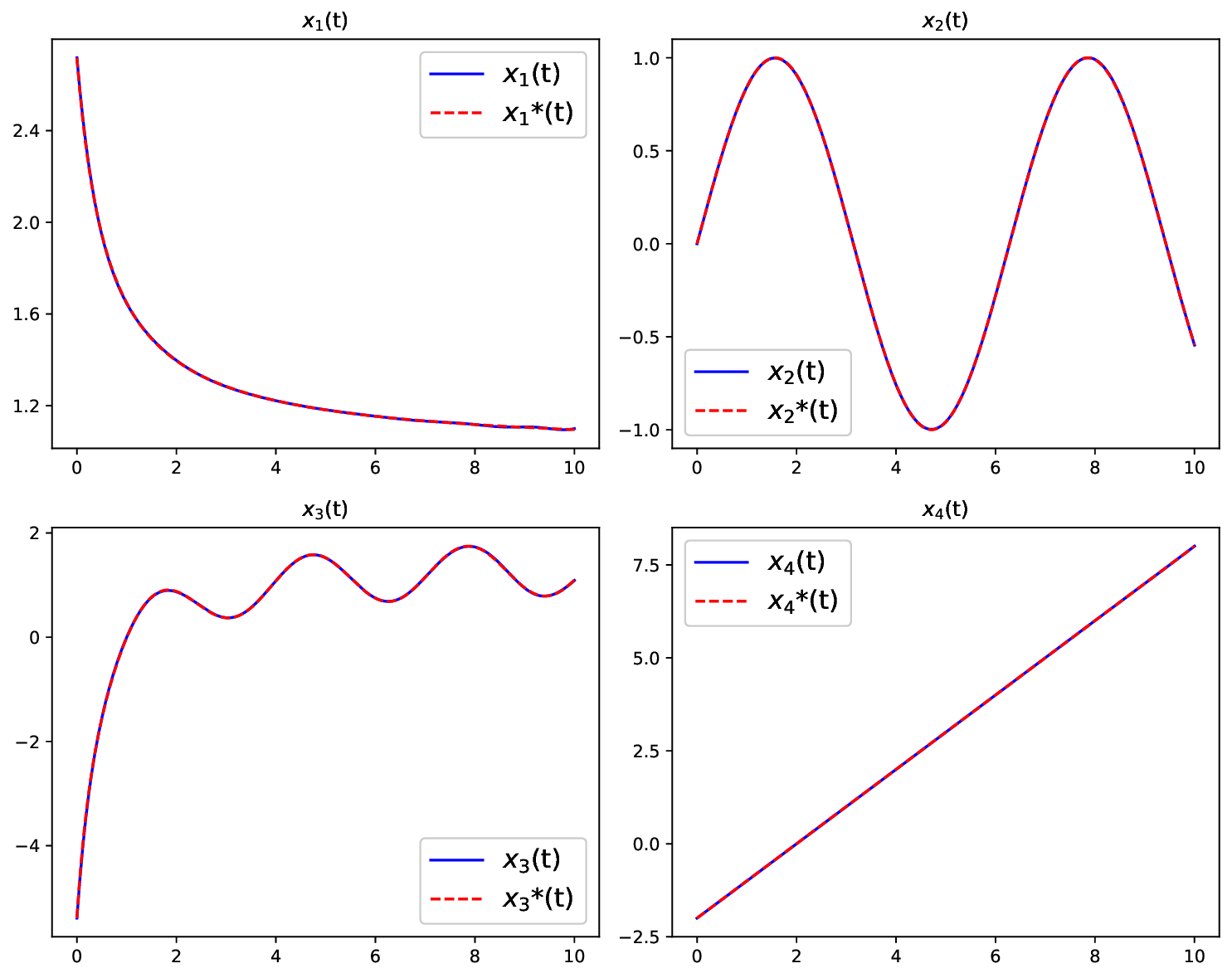}}
        \centerline{A}
	\end{minipage}
\hspace{0.2cm}
    \begin{minipage}{0.46\linewidth}
		\vspace{3pt}
		\centerline{\includegraphics[width=\textwidth]{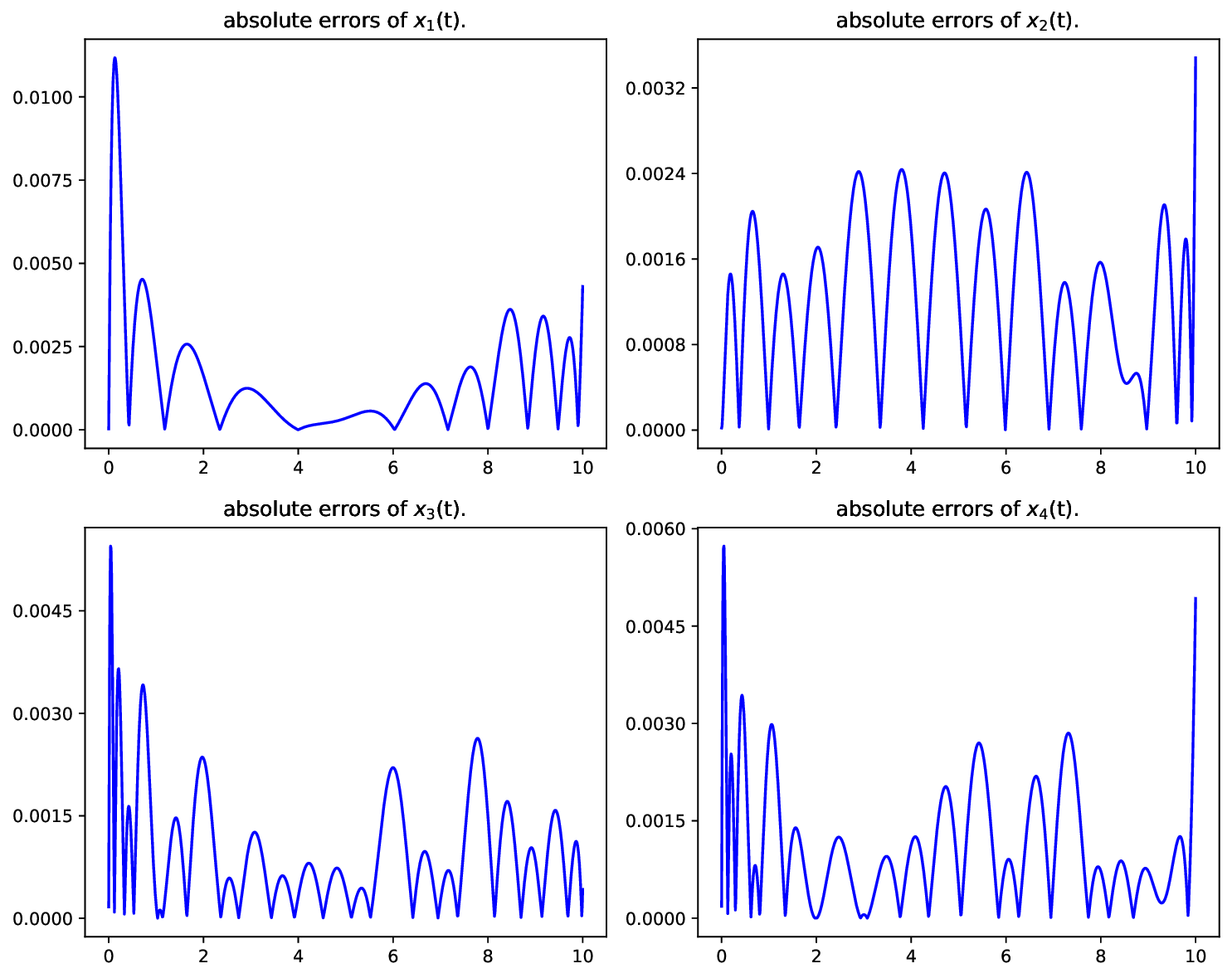}}
        \centerline{B}
	\end{minipage}
    \caption{\small (Color online) System (\ref{example-8}): (A). Exact solutions and predicted solutions; (B). Absolute errors of predicted solutions.}
    \label{fig:mte10}
\end{figure}

Next, we employ the idea of HANN-1 to solve system (\ref{example-8}). Specifically, we first sample 1000 points from $t\in[0,10]$ by the Latin hypercube sampling as the input data, and then deploy the L-BFGS optimizer to minimize the loss function $\mathcal{L}(\Theta) = \sum_{j=1}^4\left[\mathcal{L}_j(\Theta) + \mathcal{L}_{f_j}(\Theta)\right]$, where
\begin{eqnarray} \label{eq:MSE-time}
\no &&\mathcal{L}_j(\Theta)= \left|\hat{x}_j(\Theta; 0) - x_j^*(0) \right|^2, ~~~\mathcal{L}_{f_j}(\Theta)= \frac{1}{1000}\sum_{k=1}^{1000}\left|f_j(\hat{x}(\Theta; t_k),t_k) \right|^2,
\end{eqnarray}
and $f_j$ is given by system (\ref{example-8}). Note that the true value $ x_j^*(0)$ is obtained by solving system (\ref{example-8}) enforced by $t=0$, not from the true solutions (\ref{true-so}) since they are unknown in real scenario. Figure \ref{fig:mte10}(A) shows that the dotted red curves of predicted solutions fits the solid blue curves of true solutions perfectly for each $x_j\,(j=1,2,3,4)$. Moreover, in order to further verify the accuracies of predicted solutions, Figure \ref{fig:mte10}(B) displays their absolute errors by means of true solutions  (\ref{true-so}), where the absolute errors for $x_2(t),\, x_3(t)$ and $x_4(t)$ reach $10^{-3}$ order of magnitude, while the maximum occurs for $x_1(t)$ and arrives at $1.12\times10^{-2}$ which is still superior to the predicted results by the finite-time recurrent neural networks where comparatively big deviations occur near the initial zero point \cite{xzl-2019}.

\section{Further discussions}
\subsection{Systematic investigations on the six equations}
At the end of this section, we systematically evaluate the computational costs and the stabilizations of HANN-1 for the six benchmark equations, where convergent-control parameter $\gamma=0.01$. Specifically, Table \ref{tab:time costs} shows the computational costs of HANN-1 for the six equations, less than one second for the single equation and two transcendental systems, and more than thirty seconds for the time-varying system (\ref{example-8}) since it needs to compute $x_i(t)\,(i=1,\dots,4)$ in larger domain $t\in[0,10]$.  
\begin{table}[ht]
 \captionsetup{width=.45\textwidth}
\centering
\small
\caption{\small Computational costs of HANN-1 for the six examples. The data (means  $\pm$  standard error) are the results of three random trials. }
\label{tab:time costs}
\renewcommand{\arraystretch}{1.1}
\begin{tabular}{cc@{}}
\toprule
Example & Computational cost (s)\\ \midrule
 Equation (\ref{exa-1})  &	0.40  $\pm$  0.012\\
System (\ref{example-2})&1.57  $\pm$  0.023\\
System (\ref{example-3}) &	0.43  $\pm$  0.033\\
System (\ref{example-6})&	3.45  $\pm$  0.43\\
System (\ref{example-7})&1.37  $\pm$  0.12\\
System (\ref{example-8}) &	39.07  $\pm$  5.63\\
\bottomrule
\end{tabular}
\end{table}

Furthermore, we investigate the effects of HANN-1 for the six examples by considering different numbers of collocation points and neurons in each layer. Consequently, the experimental results in both Table \ref{tab:nf} and Table \ref{tab:l-n} show that proper selections of collocation point numbers and network architectures are crucial to the good performances of HANN-1 for learning solutions of the six examples. Specifically, the best results occur at 50 and 500 collocation points in Table \ref{tab:nf} and two and four hidden layers with 10, 20 and 40 neurons per layer in Table \ref{tab:l-n}, respectively. Moreover, only using HANN-1, the equation residuals arrive at $10^{-3}$ order of magnitude except for system (\ref{example-6}), even reaching $10^{-4}$ in two hidden layers with 10 neurons per layer for equation (\ref{exa-1}).
\begin{table}[ht]
 \captionsetup{width=\textwidth}
\centering
\small
\footnotesize
\caption{\small Numerical results of the six benchmark equations by HANN-1, consisting of four hidden layers and forty neurons in each layer, under different numbers of collocation points, where the data (means  $\pm$  standard error) are the results of three random trials and the bold data denote the best one.}
\label{tab:nf}
\renewcommand{\arraystretch}{1.1}
\begin{tabular}{cccccc}
\toprule
\resizebox{1.5cm}{0.5cm}{\diagbox{Eqns}{\hspace{-0.1cm}\( N_f \)}} & 5 & 50&500&1000 \\ \midrule
 Equation (\ref{exa-1}) &1.62e-03 $\pm$ 1.49e-03&\textbf{7.71e-04 $\pm$ 3.49e-04}&4.16e-03 $\pm$ 5.84e-04&1.13e-03 $\pm$ 8.09e-04\\
 System (\ref{example-2})&3.92e-03 $\pm$ 1.42e-03	&3.13e-03 $\pm$ 1.84e-03	&\textbf{1.95e-03 $\pm$ 1.74e-03}	&2.03e-03 $\pm$ 1.08e-03\\
 System (\ref{example-3}) &1.64e-01 $\pm$ 1.54e-01&9.54e-03 $\pm$ 3.17e-03&\textbf{8.06e-03 $\pm$ 1.92e-03}&7.86e-03 $\pm$ 3.22e-03\\
 System (\ref{example-6})&1.04e-01 $\pm$ 7.17e-02&\textbf{2.54e-02 $\pm$ 4.79e-03}&2.93e-02 $\pm$ 6.17e-03&3.45e-02 $\pm$ 1.19e-02\\
 System (\ref{example-7})&1.48e-02 $\pm$ 2.74e-03&1.73e-02 $\pm$ 7.32e-03&\textbf{5.98e-03 $\pm$ 2.37e-03}&9.41e-03 $\pm$ 2.86e-03\\
 System (\ref{example-8})&6.73e-01 $\pm$ 5.49e-01&1.84e-03 $\pm$ 5.42e-04&\textbf{1.04e-03 $\pm$ 3.86e-04}&1.52e-03 $\pm$ 3.57e-04\\
\bottomrule
\end{tabular}
\end{table}
\begin{table}[ht]
 \captionsetup{width=\textwidth}
\centering
\footnotesize
\caption{\small Numerical results of the six benchmark equations by HANN-1 with 500 collocation points under different network architecture, where the bold data denote the best one. The six data (means  $\pm$  standard error) in each pair of layer and neuron are the results of three random trials and for equation (\ref{exa-1}) and systems (\ref{example-2}), (\ref{example-3}), (\ref{example-6}), (\ref{example-7}) and (\ref{example-8}) in turn.}
\label{tab:l-n}
\renewcommand{\arraystretch}{1.1}
\begin{tabular}{ccccccc}
\toprule
\resizebox{1.5cm}{0.5cm}{\diagbox{Layer}{\hspace{-0.1cm}Neuron}}& 10 & 20&40&80 \\ \midrule
&\textbf{8.34e-04 $\pm$ 8.37e-05}&	2.09e-03 $\pm$ 1.44e-03&1.52e-03 $\pm$ 2.25e-04	&1.96e-03 $\pm$ 4.17e-04\\
&5.74e-03 $\pm$ 8.36e-04&	3.72e-03 $\pm$ 3.12e-03	&1.78e-03 $\pm$ 3.30e-04	&2.73e-03 $\pm$ 2.25e-03\\
&1.01e-02 $\pm$ 5.02e-03	&8.99e-03 $\pm$ 2.41e-03	&\textbf{1.02e-03 $\pm$ 8.24e-03}	&2.29e-03 $\pm$ 4.48e-04\\
2&2.78e-02 $\pm$ 4.48e-03	&\textbf{3.35e-03 $\pm$ 6.56e-04}	&4.93e-02 $\pm$ 2.45e-03	&4.89e-02 $\pm$ 1.16e-02\\
&1.45e-02 $\pm$ 2.71e-03	&9.95e-03 $\pm$ 3.39e-03	&1.19e-02 $\pm$ 4.36e-03	&4.87e-03 $\pm$ 1.77e-03\\
&2.43e-03 $\pm$ 1.69e-03	&2.08e-03 $\pm$ 4.27e-04	&2.12e-03 $\pm$ 5.52e-04	&1.52e-03 $\pm$ 3.57e-04\\
\hline
&1.78e-03 $\pm$ 1.12e-03	&1.57e-03 $\pm$ 1.19e-03	&1.86e-03 $\pm$ 1.29e-03	&9.52e-04 $\pm$ 6.46e-04\\
&3.30e-03 $\pm$ 1.67e-03	&1.45e-03 $\pm$ 9.79e-04	&2.29e-03 $\pm$ 1.09e-03	&1.11e-03 $\pm$ 2.75e-04\\
&5.11e-03 $\pm$ 2.31e-03	&2.86e-03 $\pm$ 1.80e-03	&9.53e-03 $\pm$ 1.31e-03	&4.54e-03 $\pm$ 2.53e-03\\
4&1.84e-02 $\pm$ 6.07e-03	&4.37e-02 $\pm$ 7.03e-03	&3.36e-02 $\pm$ 1.29e-03	&3.63e-02 $\pm$ 2.42e-03\\
&1.44e-02 $\pm$ 2.30e-03	&1.11e-02 $\pm$ 2.41e-03	&\textbf{6.21e-03 $\pm$ 2.64e-03}	&1.01e-02 $\pm$ 2.94e-03\\
&\textbf{1.12e-03 $\pm$ 5.52e-04}	&1.72e-03 $\pm$ 6.22e-04	&3.48e-03 $\pm$ 6.05e-04	&1.94e-03 $\pm$ 3.83e-05\\
\hline
&9.14e-04 $\pm$ 7.06e-04	&2.16e-03 $\pm$ 1.28e-03	&2.15e-03 $\pm$ 1.39e-03	&1.32e-03 $\pm$ 1.28e-03\\
&2.88e-03 $\pm$ 1.99e-03	&1.87e-03 $\pm$ 5.03e-04	&\textbf{9.62e-04 $\pm$ 2.85e-04}	&4.85e-03 $\pm$ 2.61e-03\\
&3.64e-03 $\pm$ 1.61e-03	&4.31e-03 $\pm$ 9.12e-04	&2.76e-03 $\pm$ 1.49e-03	&9.72e-03 $\pm$ 7.91e-03\\
6&1.43e-02 $\pm$ 1.49e-03	&1.35e-02 $\pm$ 4.83e-03	&2.19e-02 $\pm$ 1.62e-03	&2.79e-02 $\pm$ 7.44e-03\\
&1.34e-02 $\pm$ 1.74e-03	&1.54e-02 $\pm$ 4.43e-03	&9.75e-03 $\pm$ 5.57e-03	&1.52e-02 $\pm$ 7.52e-03\\
&2.01e-03 $\pm$ 1.19e-03	&1.22e-03 $\pm$ 3.54e-05	&1.45e-03 $\pm$ 1.12e-05	&1.46e-03 $\pm$ 3.09e-04\\
\bottomrule
\end{tabular}
\end{table}
\subsection{Conclusion}
We introduce the HANN method for solving the nonlinear algebraic equations and validate the efficacy of HANN by performing numerical experiments for the six benchmark equations which have been studied many times in literatures. It is worthy of saying that the HANN-1 introduces the homotopy parameter $t\in[0,1]$ as the input and direct output the final solutions of equation which is completely different with the known neural network methods for solving algebraic equations. Moreover, the HANN-2 iteratively employs HANN-1 to refine the learned solutions and the first HANN-1 training phase dominates the total computational time, shown in the loss history of HANN-2 in Figure \ref{fig:sideex-1122}(D), thus the computational costs by HANN-2 will not increase significantly. Therefore, in practice, we recommend first using HANN-1 with sufficient number of random initial values to rapidly identify potential solution candidates, then leverage these approximate solutions as the initial values of HANN-2 or the Python's built-in Fsolve function to generate high-accuracy solutions.

\section*{Acknowledgements}
The paper is supported by the Beijing Natural Science Foundation (No. 1222014) and the National Natural Science Foundation of China (No. 11671014).
\vspace{0.12cm}\\ \textbf{Declarations of interest: The authors have no conflicts to disclose.}

\end{document}